\documentclass[12pt]{article}
\usepackage[amsmath]{e-jc}

\usepackage{graphicx}

\usepackage{subeqnarray, eepicemu, url,cite,bm}
\usepackage[T1]{fontenc}

\dateline{Jun 20, 2023}{Feb 29, 2024}{TBD}


\Copyright{The authors. Released under the CC BY-ND license (International 4.0).}

\title{A remark on continued fractions for permutations \\
         and D-permutations with a weight $-1$ per cycle}

\author{
Bishal Deb\authornote{1,2} \and
Alan D. Sokal\authornote{1,3}
}

\authortext{1}{
Department of Mathematics, University College London,
London WC1E 6BT, UK
(\email{bishal@gonitsora.com}, \email{sokal@nyu.edu}).
}

\authortext{2}{
Sorbonne Universit\'e and Universit\'e Paris Cit\'e, CNRS,\\ \hspace{5mm} Laboratoire de Probabilit\'es,
Statistique et Mod\'elisation, Paris, France.
}

\authortext{3}{Department of Physics, New York University, New York, NY 10003, USA.
}

\begin{document}

\maketitle


\begin{abstract}
We show that very simple continued fractions can be obtained
for the ordinary generating functions enumerating
permutations or D-permutations with a large number of independent statistics,
when each cycle is given a weight $-1$.
The proof is based on a simple lemma
relating the number of cycles modulo 2
to the numbers of fixed points, cycle peaks (or cycle valleys), and crossings.
\end{abstract}

\renewcommand{\theenumi}{\alph{enumi}}
\renewcommand{\labelenumi}{(\theenumi)}
\def\eop{\hbox{\kern1pt\vrule height6pt width4pt
depth1pt\kern1pt}\medskip}
\def\prf{\par\noindent{\bf Proof.\enspace}\rm}
\def\rmk{\par\medskip\noindent{\bf Remark\enspace}\rm}

\newcommand{\textbfit}[1]{\textbf{\textit{#1}}}

\newcommand{\bigdash}{%
\smallskip\begin{center} \rule{5cm}{0.1mm} \end{center}\smallskip}

\newcommand{\safepar}{ {\protect\hfill\protect\break\hspace*{5mm}} }

\newcommand{\be}{\begin{equation}}
\newcommand{\ee}{\end{equation}}
\newcommand{\<}{\langle}
\renewcommand{\>}{\rangle}
\newcommand{\widebar}{\overline}
\def\reff#1{(\protect\ref{#1})}
\def\spose#1{\hbox to 0pt{#1\hss}}
\def\ltapprox{\mathrel{\spose{\lower 3pt\hbox{$\mathchar"218$}}
    \raise 2.0pt\hbox{$\mathchar"13C$}}}
\def\gtapprox{\mathrel{\spose{\lower 3pt\hbox{$\mathchar"218$}}
    \raise 2.0pt\hbox{$\mathchar"13E$}}}
\def\textprime{${}^\prime$}
\def\proof{\par\medskip\noindent{\sc Proof.\ }}
\def\firstproof{\par\medskip\noindent{\sc First Proof.\ }}
\def\secondproof{\par\medskip\noindent{\sc Second Proof.\ }}
\def\alternateproof{\par\medskip\noindent{\sc Alternate Proof.\ }}
\def\algebraicproof{\par\medskip\noindent{\sc Algebraic Proof.\ }}
\def\combinatorialproof{\par\medskip\noindent{\sc Combinatorial Proof.\ }}
\def\topologicalproof{\par\medskip\noindent{\sc Topological Proof.\ }}
\def\proofof#1{\bigskip\noindent{\sc Proof of #1.\ }}
\def\firstproofof#1{\bigskip\noindent{\sc First Proof of #1.\ }}
\def\secondproofof#1{\bigskip\noindent{\sc Second Proof of #1.\ }}
\def\thirdproofof#1{\bigskip\noindent{\sc Third Proof of #1.\ }}
\def\algebraicproofof#1{\bigskip\noindent{\sc Algebraic Proof of #1.\ }}
\def\combinatorialproofof#1{\bigskip\noindent{\sc Combinatorial Proof of #1.\ }}
\def\sketchofproof{\par\medskip\noindent{\sc Sketch of proof.\ }}
\renewcommand{\qed}{ $\square$ \bigskip}
\newcommand{\myendremark}{ $\blacksquare$ \bigskip}
\def\half{ {1 \over 2} }
\def\third{ {1 \over 3} }
\def\twothird{ {2 \over 3} }
\def\smfrac#1#2{{\textstyle{#1\over #2}}}
\def\smhalf{ {\smfrac{1}{2}} }
\newcommand{\real}{\mathop{\rm Re}\nolimits}
\renewcommand{\Re}{\mathop{\rm Re}\nolimits}
\newcommand{\imag}{\mathop{\rm Im}\nolimits}
\renewcommand{\Im}{\mathop{\rm Im}\nolimits}
\newcommand{\sgn}{\mathop{\rm sgn}\nolimits}
\newcommand{\tr}{\mathop{\rm tr}\nolimits}
\newcommand{\supp}{\mathop{\rm supp}\nolimits}
\newcommand{\disc}{\mathop{\rm disc}\nolimits}
\newcommand{\diag}{\mathop{\rm diag}\nolimits}
\newcommand{\tridiag}{\mathop{\rm tridiag}\nolimits}
\newcommand{\AZ}{\mathop{\rm AZ}\nolimits}
\newcommand{\NC}{\mathop{\rm NC}\nolimits}
\newcommand{\PF}{{\rm PF}}
\newcommand{\rk}{\mathop{\rm rk}\nolimits}
\newcommand{\perm}{\mathop{\rm perm}\nolimits}
\def\hboxscript#1{ {\hbox{\scriptsize\em #1}} }
\renewcommand{\emptyset}{\varnothing}
\newcommand{\eqdef}{\stackrel{\rm def}{=}}

\newcommand{\restrict}{\upharpoonright}

\newcommand{\compinv}{{\langle -1 \rangle}}   

\newcommand{\scra}{{\mathcal{A}}}
\newcommand{\scrb}{{\mathcal{B}}}
\newcommand{\scrc}{{\mathcal{C}}}
\newcommand{\bfscra}{{\bm{\mathcal{A}}}}
\newcommand{\bfscrb}{{\bm{\mathcal{B}}}}
\newcommand{\bfscrc}{{\bm{\mathcal{C}}}}
\newcommand{\bfscrap}{{\bm{\mathcal{A}'}}}
\newcommand{\bfscrbp}{{\bm{\mathcal{B}'}}}
\newcommand{\bfscrcp}{{\bm{\mathcal{C}'}}}
\newcommand{\bfscrapp}{{\bm{\mathcal{A}''}}}
\newcommand{\bfscrbpp}{{\bm{\mathcal{B}''}}}
\newcommand{\bfscrcpp}{{\bm{\mathcal{C}''}}}
\newcommand{\scrd}{{\mathcal{D}}}
\newcommand{\scre}{{\mathcal{E}}}
\newcommand{\scrf}{{\mathcal{F}}}
\newcommand{\scrg}{{\mathcal{G}}}
\newcommand{\scrgg}{{\mathscr{g}}}  
\newcommand{\scrh}{{\mathcal{H}}}
\newcommand{\scri}{{\mathcal{I}}}
\newcommand{\scrj}{{\mathcal{J}}}
\newcommand{\scrk}{{\mathcal{K}}}
\newcommand{\scrl}{{\mathcal{L}}}
\newcommand{\scrm}{{\mathcal{M}}}
\newcommand{\scrn}{{\mathcal{N}}}
\newcommand{\scro}{{\mathcal{O}}}
\newcommand\scroo{
  \mathchoice
    {{\scriptstyle\mathcal{O}}}
    {{\scriptstyle\mathcal{O}}}
    {{\scriptscriptstyle\mathcal{O}}}
    {\scalebox{0.6}{$\scriptscriptstyle\mathcal{O}$}}
  }
\newcommand{\scrp}{{\mathcal{P}}}
\newcommand{\scrq}{{\mathcal{Q}}}
\newcommand{\scrr}{{\mathcal{R}}}
\newcommand{\scrs}{{\mathcal{S}}}
\newcommand{\scrss}{{\mathscr{s}}}  
\newcommand{\scrt}{{\mathcal{T}}}
\newcommand{\scrv}{{\mathcal{V}}}
\newcommand{\scrw}{{\mathcal{W}}}
\newcommand{\scrz}{{\mathcal{Z}}}
\newcommand{\SP}{{\mathcal{SP}}}
\newcommand{\ST}{{\mathcal{ST}}}

\newcommand{\bfa}{{\mathbf{a}}}
\newcommand{\bfb}{{\mathbf{b}}}
\newcommand{\bfc}{{\mathbf{c}}}
\newcommand{\bfd}{{\mathbf{d}}}
\newcommand{\bfe}{{\mathbf{e}}}
\newcommand{\bfh}{{\mathbf{h}}}
\newcommand{\bfj}{{\mathbf{j}}}
\newcommand{\bfi}{{\mathbf{i}}}
\newcommand{\bfk}{{\mathbf{k}}}
\newcommand{\bfl}{{\mathbf{l}}}
\newcommand{\bfL}{{\mathbf{L}}}
\newcommand{\bfm}{{\mathbf{m}}}
\newcommand{\bfn}{{\mathbf{n}}}
\newcommand{\bfp}{{\mathbf{p}}}
\newcommand{\bfr}{{\mathbf{r}}}
\newcommand{\bfu}{{\mathbf{u}}}
\newcommand{\bfv}{{\mathbf{v}}}
\newcommand{\bfw}{{\mathbf{w}}}
\newcommand{\bfx}{{\mathbf{x}}}
\newcommand{\bfy}{{\mathbf{y}}}
\newcommand{\bfz}{{\mathbf{z}}}
\renewcommand{\k}{{\mathbf{k}}}
\newcommand{\n}{{\mathbf{n}}}
\newcommand{\vv}{{\mathbf{v}}}
\newcommand{\bv}{{\mathbf{v}}}
\newcommand{\w}{{\mathbf{w}}}
\newcommand{\x}{{\mathbf{x}}}
\newcommand{\y}{{\mathbf{y}}}
\newcommand{\cc}{{\mathbf{c}}}
\newcommand{\zero}{{\mathbf{0}}}
\newcommand{\one}{{\mathbf{1}}}
\newcommand{\bmm}{{\mathbf{m}}}

\newcommand{\ahat}{{\widehat{a}}}
\newcommand{\vhat}{{\widehat{v}}}
\newcommand{\yhat}{{\widehat{y}}}
\newcommand{\phat}{{\widehat{p}}}
\newcommand{\qhat}{{\widehat{q}}}
\newcommand{\Zhat}{{\widehat{Z}}}
\newcommand{\myhat}{{\widehat{\;}}}
\newcommand{\vtilde}{{\widetilde{v}}}
\newcommand{\ytilde}{{\widetilde{y}}}

\newcommand{\C}{{\mathbb C}}
\newcommand{\D}{{\mathbb D}}
\newcommand{\Z}{{\mathbb Z}}
\newcommand{\N}{{\mathbb N}}
\newcommand{\Q}{{\mathbb Q}}
\newcommand{\PP}{{\mathbb P}}
\newcommand{\R}{{\mathbb R}}
\newcommand{\RR}{{\mathbb R}}
\newcommand{\E}{{\mathbb E}}

\newcommand{\Sym}{{\mathfrak{S}}}
\newcommand{\SymB}{{\mathfrak{B}}}
\newcommand{\Cyc}{{\mathfrak{C}}}
\newcommand{\Altcyc}{{\scra\mathfrak{C}}}
\newcommand{\Alt}{{\mathrm{Alt}}}
\newcommand{\dperm}{{\mathfrak{D}}}
\newcommand{\dcycle}{{\mathfrak{DC}}}

\newcommand{\germanA}{{\mathfrak{A}}}
\newcommand{\germanB}{{\mathfrak{B}}}
\newcommand{\germanQ}{{\mathfrak{Q}}}
\newcommand{\germanh}{{\mathfrak{h}}}

\newcommand{\myle}{\preceq}
\newcommand{\myge}{\succeq}
\newcommand{\mygt}{\succ}

\newcommand{\B}{{\sf B}}
\newcommand{\OB}{B^{\rm ord}}
\newcommand{\OS}{{\sf OS}}
\newcommand{\OO}{{\sf O}}
\newcommand{\OSP}{{\sf OSP}}
\newcommand{\Eu}{{\sf Eu}}
\newcommand{\ERR}{{\sf ERR}}
\newcommand{\sfB}{{\sf B}}
\newcommand{\sfD}{{\sf D}}
\newcommand{\sfE}{{\sf E}}
\newcommand{\sfG}{{\sf G}}
\newcommand{\sfJ}{{\sf J}}
\newcommand{\sfL}{{\sf L}}
\newcommand{\sfLhat}{{\widehat{{\sf L}}}}
\newcommand{\sfLcheck}{{\widecheck{{\sf L}}}}
\newcommand{\sfLtilde}{{\widetilde{{\sf L}}}}
\newcommand{\sfP}{{\sf P}}
\newcommand{\sfQ}{{\sf Q}}
\newcommand{\sfS}{{\sf S}}
\newcommand{\sfT}{{\sf T}}
\newcommand{\sfW}{{\sf W}}
\newcommand{\sfMV}{{\sf MV}}
\newcommand{\AMV}{{\sf AMV}}
\newcommand{\BM}{{\sf BM}}
\newcommand{\emIB}{B^{\rm irr}}
\newcommand{\emIP}{P^{\rm irr}}
\newcommand{\emOB}{B^{\rm ord}}
\newcommand{\emCB}{B^{\rm cyc}}
\newcommand{\emSC}{P^{\rm cyc}}

\newcommand{\lev}{{\rm lev}}
\newcommand{\stat}{{\rm stat}}
\newcommand{\cyc}{{\rm cyc}}
\newcommand{\mysteryone}{{\rm mys1}}
\newcommand{\mysterytwo}{{\rm mys2}}
\newcommand{\Asc}{{\rm Asc}}
\newcommand{\asc}{{\rm asc}}
\newcommand{\Des}{{\rm Des}}
\newcommand{\des}{{\rm des}}
\newcommand{\Exc}{{\rm Exc}}

\newcommand{\EArec}{{\rm EArec}}
\newcommand{\earec}{{\rm earec}}
\newcommand{\recarec}{{\rm recarec}}
\newcommand{\erec}{{\rm erec}}
\newcommand{\nonrec}{{\rm nonrec}}
\newcommand{\nrar}{{\rm nrar}}
\newcommand{\Ereccval}{{\rm Ereccval}}
\newcommand{\ereccval}{{\rm ereccval}}
\newcommand{\ereccvalodd}{{\rm ereccvalodd}}
\newcommand{\ereccvaleven}{{\rm ereccvaleven}}
\newcommand{\ereccdrise}{{\rm ereccdrise}}
\newcommand{\Eareccpeak}{{\rm Eareccpeak}}
\newcommand{\eareccpeak}{{\rm eareccpeak}}
\newcommand{\Eareccpeakodd}{{\rm Eareccpeakodd}}
\newcommand{\Eareccpeakeven}{{\rm Eareccpeakeven}}
\newcommand{\eareccpeakodd}{{\rm eareccpeakodd}}
\newcommand{\eareccpeakeven}{{\rm eareccpeakeven}}
\newcommand{\eareccdfall}{{\rm eareccdfall}}
\newcommand{\eareccval}{{\rm eareccval}}
\newcommand{\ereccpeak}{{\rm ereccpeak}}
\newcommand{\Nrpeak}{{\rm Nrpeak}}
\newcommand{\rar}{{\rm rar}}
\newcommand{\evenrar}{{\rm evenrar}}
\newcommand{\oddrar}{{\rm oddrar}}
\newcommand{\nrcpeak}{{\rm nrcpeak}}
\newcommand{\nrcpeakodd}{{\rm nrcpeakodd}}
\newcommand{\nrcpeakeven}{{\rm nrcpeakeven}}
\newcommand{\nrcval}{{\rm nrcval}}
\newcommand{\nrcvalodd}{{\rm nrcvalodd}}
\newcommand{\nrcvaleven}{{\rm nrcvaleven}}
\newcommand{\nrcdrise}{{\rm nrcdrise}}
\newcommand{\nrcdfall}{{\rm nrcdfall}}
\newcommand{\nrfix}{{\rm nrfix}}
\newcommand{\Evenfix}{{\rm Evenfix}}
\newcommand{\evenfix}{{\rm evenfix}}
\newcommand{\Oddfix}{{\rm Oddfix}}
\newcommand{\oddfix}{{\rm oddfix}}
\newcommand{\evennrfix}{{\rm evennrfix}}
\newcommand{\oddnrfix}{{\rm oddnrfix}}
\newcommand{\Cpeak}{{\rm Cpeak}}
\newcommand{\cpeak}{{\rm cpeak}}
\newcommand{\Cval}{{\rm Cval}}
\newcommand{\cval}{{\rm cval}}
\newcommand{\cvaleven}{{\rm cvaleven}}
\newcommand{\cvalodd}{{\rm cvalodd}}
\newcommand{\Cdasc}{{\rm Cdasc}}
\newcommand{\cdasc}{{\rm cdasc}}
\newcommand{\Cddes}{{\rm Cddes}}
\newcommand{\cddes}{{\rm cddes}}
\newcommand{\Cdrise}{{\rm Cdrise}}
\newcommand{\cdrise}{{\rm cdrise}}
\newcommand{\Cdfall}{{\rm Cdfall}}
\newcommand{\cdfall}{{\rm cdfall}}

\newcommand{\maxpeak}{{\rm maxpeak}}
\newcommand{\nmaxpeak}{{\rm nmaxpeak}}
\newcommand{\minval}{{\rm minval}}
\newcommand{\nminval}{{\rm nminval}}

\newcommand{\exc}{{\rm exc}}
\newcommand{\excee}{{\rm excee}}
\newcommand{\exceo}{{\rm exceo}}
\newcommand{\excoe}{{\rm excoe}}
\newcommand{\excoo}{{\rm excoo}}
\newcommand{\erecexcoe}{{\rm erecexcoe}}
\newcommand{\erecexcoo}{{\rm erecexcoo}}
\newcommand{\nrexcoe}{{\rm nrexcoe}}
\newcommand{\nrexcoo}{{\rm nrexcoo}}
\newcommand{\aexc}{{\rm aexc}}
\newcommand{\aexcee}{{\rm aexcee}}
\newcommand{\aexceo}{{\rm aexceo}}
\newcommand{\aexcoe}{{\rm aexcoe}}
\newcommand{\aexcoo}{{\rm aexcoo}}
\newcommand{\earecaexcee}{{\rm earecaexcee}}
\newcommand{\earecaexceo}{{\rm earecaexceo}}
\newcommand{\nraexcee}{{\rm nraexcee}}
\newcommand{\nraexceo}{{\rm nraexceo}}
\newcommand{\Fix}{{\rm Fix}}
\newcommand{\fix}{{\rm fix}}
\newcommand{\fixe}{{\rm fixe}}
\newcommand{\fixo}{{\rm fixo}}
\newcommand{\bfix}{{\mathbf{fix}}}
\newcommand{\rare}{{\rm rare}}
\newcommand{\raro}{{\rm raro}}
\newcommand{\nrfixe}{{\rm nrfixe}}
\newcommand{\nrfixo}{{\rm nrfixo}}
\newcommand{\xee}{x_{\rm ee}}
\newcommand{\xeo}{x_{\rm eo}}
\newcommand{\uee}{u_{\rm ee}}
\newcommand{\ueo}{u_{\rm eo}}
\newcommand{\yoo}{y_{\rm oo}}
\newcommand{\yoe}{y_{\rm oe}}
\newcommand{\voo}{v_{\rm oo}}
\newcommand{\voe}{v_{\rm oe}}
\newcommand{\zo}{z_{\rm o}}
\newcommand{\ze}{z_{\rm e}}
\newcommand{\wo}{w_{\rm o}}
\newcommand{\we}{w_{\rm e}}
\newcommand{\so}{s_{\rm o}}
\newcommand{\se}{s_{\rm e}}

\newcommand{\xo}{x_{\rm o}}
\newcommand{\xe}{x_{\rm e}}
\newcommand{\yo}{y_{\rm o}}
\newcommand{\ye}{y_{\rm e}}
\newcommand{\uo}{u_{\rm o}}
\newcommand{\ue}{u_{\rm e}}
\newcommand{\vo}{v_{\rm o}}
\newcommand{\ve}{v_{\rm e}}

\newcommand{\Wex}{{\rm Wex}}
\newcommand{\wex}{{\rm wex}}
\newcommand{\lrmax}{{\rm lrmax}}
\newcommand{\rlmax}{{\rm rlmax}}
\newcommand{\Rec}{{\rm Rec}}
\newcommand{\rec}{{\rm rec}}
\newcommand{\Arec}{{\rm Arec}}
\newcommand{\Arecpeak}{{\rm Arecpeak}}
\newcommand{\arec}{{\rm arec}}
\newcommand{\Even}{{\rm Even}}
\newcommand{\Odd}{{\rm Odd}}
\newcommand{\ERec}{{\rm ERec}}
\newcommand{\Val}{{\rm Val}}
\newcommand{\Peak}{{\rm Peak}}
\newcommand{\dasc}{{\rm dasc}}
\newcommand{\ddes}{{\rm ddes}}
\newcommand{\inv}{{\rm inv}}
\newcommand{\maj}{{\rm maj}}
\newcommand{\rs}{{\rm rs}}
\newcommand{\cross}{{\rm cr}}
\newcommand{\crosshat}{{\widehat{\rm cr}}}
\newcommand{\nest}{{\rm ne}}
\newcommand{\ucross}{{\rm ucross}}
\newcommand{\Ucross}{{\rm Ucross}}
\newcommand{\ucrosscval}{{\rm ucrosscval}}
\newcommand{\Ucrosscval}{{\rm Ucrosscval}}
\newcommand{\ucrosscpeak}{{\rm ucrosscpeak}}
\newcommand{\Ucrosscpeak}{{\rm Ucrosscpeak}}
\newcommand{\ucrosscdrise}{{\rm ucrosscdrise}}
\newcommand{\Ucrosscdrise}{{\rm Ucrosscdrise}}
\newcommand{\lcross}{{\rm lcross}}
\newcommand{\Lcross}{{\rm Lcross}}
\newcommand{\lcrosscpeak}{{\rm lcrosscpeak}}
\newcommand{\Lcrosscpeak}{{\rm Lcrosscpeak}}
\newcommand{\lcrosscval}{{\rm lcrosscval}}
\newcommand{\Lcrosscval}{{\rm Lcrosscval}}
\newcommand{\lcrosscdfall}{{\rm lcrosscdfall}}
\newcommand{\Lcrosscdfall}{{\rm Lcrosscdfall}}
\newcommand{\unest}{{\rm unest}}
\newcommand{\Unest}{{\rm Unest}}
\newcommand{\unestcval}{{\rm unestcval}}
\newcommand{\Unestcval}{{\rm Unestcval}}
\newcommand{\unestcpeak}{{\rm unestcpeak}}
\newcommand{\Unestcpeak}{{\rm Unestcpeak}}
\newcommand{\unestcdrise}{{\rm unestcdrise}}
\newcommand{\Unestcdrise}{{\rm Unestcdrise}}
\newcommand{\lnest}{{\rm lnest}}
\newcommand{\Lnest}{{\rm Lnest}}
\newcommand{\lnestcpeak}{{\rm lnestcpeak}}
\newcommand{\Lnestcpeak}{{\rm Lnestcpeak}}
\newcommand{\lnestcval}{{\rm lnestcval}}
\newcommand{\Lnestcval}{{\rm Lnestcval}}
\newcommand{\lnestcdfall}{{\rm lnestcdfall}}
\newcommand{\Lnestcdfall}{{\rm Lnestcdfall}}
\newcommand{\ulev}{{\rm ulev}}
\newcommand{\llev}{{\rm llev}}
\newcommand{\ujoin}{{\rm ujoin}}
\newcommand{\ljoin}{{\rm ljoin}}
\newcommand{\psnest}{{\rm psnest}}
\newcommand{\upsnest}{{\rm upsnest}}
\newcommand{\Upsnest}{{\rm Upsnest}}
\newcommand{\lpsnest}{{\rm lpsnest}}
\newcommand{\Lpsnest}{{\rm Lpsnest}}
\newcommand{\epsnest}{{\rm epsnest}}
\newcommand{\opsnest}{{\rm opsnest}}
\newcommand{\rodd}{{\rm rodd}}
\newcommand{\reven}{{\rm reven}}
\newcommand{\lodd}{{\rm lodd}}
\newcommand{\leven}{{\rm leven}}
\newcommand{\sg}{{\rm sg}}
\newcommand{\bl}{{\rm bl}}
\newcommand{\tran}{{\rm tr}}
\newcommand{\area}{{\rm area}}
\newcommand{\ret}{{\rm ret}}
\newcommand{\peaks}{{\rm peaks}}
\newcommand{\hl}{{\rm hl}}
\newcommand{\sll}{{\rm sl}}
\newcommand{\negg}{{\rm neg}}
\newcommand{\imp}{{\rm imp}}
\newcommand{\osg}{{\rm osg}}
\newcommand{\ons}{{\rm ons}}
\newcommand{\isg}{{\rm isg}}
\newcommand{\ins}{{\rm ins}}
\newcommand{\LL}{{\rm LL}}
\newcommand{\height}{{\rm ht}}
\newcommand{\as}{{\rm as}}

\newcommand{\ba}{{\bm{a}}}
\newcommand{\bahat}{{\widehat{\bm{a}}}}
\newcommand{\bb}{{\bm{b}}}
\newcommand{\bc}{{\bm{c}}}
\newcommand{\bchat}{{\widehat{\bm{c}}}}
\newcommand{\bd}{{\bm{d}}}
\newcommand{\bee}{{\bm{e}}}
\newcommand{\beh}{{\bm{eh}}}
\newcommand{\bff}{{\bm{f}}}
\newcommand{\bg}{{\bm{g}}}
\newcommand{\bh}{{\bm{h}}}
\newcommand{\bll}{{\bm{\ell}}}
\newcommand{\bp}{{\bm{p}}}
\newcommand{\br}{{\bm{r}}}
\newcommand{\bs}{{\bm{s}}}
\newcommand{\bu}{{\bm{u}}}
\newcommand{\bw}{{\mathbf{w}}}
\newcommand{\bx}{{\bm{x}}}
\newcommand{\by}{{\bm{y}}}
\newcommand{\bz}{{\bm{z}}}
\newcommand{\bA}{{\bm{A}}}
\newcommand{\bB}{{\bm{B}}}
\newcommand{\bC}{{\bm{C}}}
\newcommand{\bE}{{\bm{E}}}
\newcommand{\bF}{{\bm{F}}}
\newcommand{\bG}{{\bm{G}}}
\newcommand{\bH}{{\bm{H}}}
\newcommand{\bI}{{\bm{I}}}
\newcommand{\bJ}{{\bm{J}}}
\newcommand{\bM}{{\bm{M}}}
\newcommand{\bN}{{\bm{N}}}
\newcommand{\bP}{{\bm{P}}}
\newcommand{\bQ}{{\bm{Q}}}
\newcommand{\bR}{{\bm{R}}}
\newcommand{\bS}{{\bm{S}}}
\newcommand{\bT}{{\bm{T}}}
\newcommand{\bW}{{\bm{W}}}
\newcommand{\bX}{{\bm{X}}}
\newcommand{\bY}{{\bm{Y}}}
\newcommand{\bIB}{{\bm{B}^{\rm irr}}}
\newcommand{\bOB}{{\bm{B}^{\rm ord}}}
\newcommand{\bOS}{{\bm{OS}}}
\newcommand{\bERR}{{\bm{ERR}}}
\newcommand{\bSP}{{\bm{SP}}}
\newcommand{\bMV}{{\bm{MV}}}
\newcommand{\bBM}{{\bm{BM}}}
\newcommand{\balpha}{{\bm{\alpha}}}
\newcommand{\balphapre}{{\bm{\alpha}^{\rm pre}}}
\newcommand{\bbeta}{{\bm{\beta}}}
\newcommand{\bgamma}{{\bm{\gamma}}}
\newcommand{\bdelta}{{\bm{\delta}}}
\newcommand{\bkappa}{{\bm{\kappa}}}
\newcommand{\bmu}{{\bm{\mu}}}
\newcommand{\bomega}{{\bm{\omega}}}
\newcommand{\bsigma}{{\bm{\sigma}}}
\newcommand{\btau}{{\bm{\tau}}}
\newcommand{\bphi}{{\bm{\phi}}}
\newcommand{\bpsi}{{\bm{\psi}}}
\newcommand{\bzeta}{{\bm{\zeta}}}
\newcommand{\bone}{{\bm{1}}}
\newcommand{\bzero}{{\bm{0}}}

\newcommand{\sfa}{{{\sf a}}}
\newcommand{\sfb}{{{\sf b}}}
\newcommand{\sfc}{{{\sf c}}}
\newcommand{\sfd}{{{\sf d}}}
\newcommand{\sfe}{{{\sf e}}}
\newcommand{\sff}{{{\sf f}}}
\newcommand{\sfg}{{{\sf g}}}
\newcommand{\sfh}{{{\sf h}}}
\newcommand{\sfi}{{{\sf i}}}
\newcommand{\bsfa}{{\mbox{\textsf{\textbf{a}}}}}
\newcommand{\bsfb}{{\mbox{\textsf{\textbf{b}}}}}
\newcommand{\bsfc}{{\mbox{\textsf{\textbf{c}}}}}
\newcommand{\bsfd}{{\mbox{\textsf{\textbf{d}}}}}
\newcommand{\bsfe}{{\mbox{\textsf{\textbf{e}}}}}
\newcommand{\bsff}{{\mbox{\textsf{\textbf{f}}}}}
\newcommand{\bsfg}{{\mbox{\textsf{\textbf{g}}}}}
\newcommand{\bsfh}{{\mbox{\textsf{\textbf{h}}}}}
\newcommand{\bsfi}{{\mbox{\textsf{\textbf{i}}}}}

\newcommand{\Cbar}{{\overline{C}}}
\newcommand{\Dbar}{{\overline{D}}}
\newcommand{\dbar}{{\overline{d}}}
\def\Ctilde{{\widetilde{C}}}
\def\Ftilde{{\widetilde{F}}}
\def\Gtilde{{\widetilde{G}}}
\def\Htilde{{\widetilde{H}}}
\def\Ptilde{{\widetilde{P}}}
\def\Chat{{\widehat{C}}}
\def\ctilde{{\widetilde{c}}}
\def\zbar{{\overline{Z}}}
\def\pitilde{{\widetilde{\pi}}}
\def\omegahat{{\widehat{\omega}}}

\newcommand{\LD}{{\mathbf{LD}}}
\newcommand{\e}{{\rm e}}
\newcommand{\ecyc}{{\rm ecyc}}
\newcommand{\epa}{{\rm epa}}
\newcommand{\iv}{{\rm iv}}
\newcommand{\pa}{{\rm pa}}
\newcommand{\pk}{{\rm p}}
\newcommand{\val}{{\rm v}}
\newcommand{\da}{{\rm da}}
\newcommand{\dd}{{\rm dd}}
\newcommand{\fp}{{\rm fp}}
\newcommand{\pkcyc}{{\rm pcyc}}
\newcommand{\valcyc}{{\rm vcyc}}
\newcommand{\dacyc}{{\rm dacyc}}
\newcommand{\ddcyc}{{\rm ddcyc}}
\newcommand{\pkpa}{{\rm ppa}}
\newcommand{\valpa}{{\rm vpa}}
\newcommand{\dapa}{{\rm dapa}}
\newcommand{\ddpa}{{\rm ddpa}}

\newcommand{\yp}{{y_\pk}}
\newcommand{\yptilde}{{\widetilde{y}_\pk}}
\newcommand{\yphat}{{\widehat{y}_\pk}}
\newcommand{\yv}{{y_\val}}
\newcommand{\yiv}{{y_\iv}}
\newcommand{\yda}{{y_\da}}
\newcommand{\ydatilde}{{\widetilde{y}_\da}}
\newcommand{\ydd}{{y_\dd}}
\newcommand{\yddtilde}{{\widetilde{y}_\dd}}
\newcommand{\yfp}{{y_\fp}}
\newcommand{\zp}{{z_\pk}}
\newcommand{\zv}{{z_\val}}
\newcommand{\zda}{{z_\da}}
\newcommand{\zdd}{{z_\dd}}

\newcommand{\sech}{{\rm sech}}

\newcommand{\sinv}{\sigma^{-1}}

%
%
\newcommand{\sn}{{\rm sn}}
\newcommand{\cn}{{\rm cn}}
\newcommand{\dn}{{\rm dn}}
\newcommand{\sm}{{\rm sm}}
\newcommand{\cm}{{\rm cm}}

%
%
\newcommand{\zfz}{ {{}_0 \! F_0} }
\newcommand{\zfo}{ {{}_0  F_1} }
\newcommand{\ofz}{ {{}_1 \! F_0} }
\newcommand{\ofo}{ {{}_1 \! F_1} }
\newcommand{\oft}{ {{}_1 \! F_2} }

%
%
\newcommand{\FHyper}[2]{ {\tensor[_{#1 \!}]{F}{_{#2}}\!} }
\newcommand{\FHYPER}[5]{ {\FHyper{#1}{#2} \!\biggl(
   \!\!\begin{array}{c} #3 \\[1mm] #4 \end{array}\! \bigg|\, #5 \! \biggr)} }
\newcommand{\tfo}{ {\FHyper{2}{1}} }
\newcommand{\tfz}{ {\FHyper{2}{0}} }
\newcommand{\threefz}{ {\FHyper{3}{0}} }
\newcommand{\FHYPERbottomzero}[3]{ {\FHyper{#1}{0} \hspace*{-0mm}\biggl(
   \!\!\begin{array}{c} #2 \\[1mm] \hbox{---} \end{array}\! \bigg|\, #3 \! \biggr)} }
\newcommand{\FHYPERtopzero}[3]{ {\FHyper{0}{#1} \hspace*{-0mm}\biggl(
   \!\!\begin{array}{c} \hbox{---} \\[1mm] #2 \end{array}\! \bigg|\, #3 \! \biggr)} }

\newcommand{\phiHyper}[2]{ {\tensor[_{#1}]{\phi}{_{#2}}} }
\newcommand{\psiHyper}[2]{ {\tensor[_{#1}]{\psi}{_{#2}}} }
\newcommand{\PhiHyper}[2]{ {\tensor[_{#1}]{\Phi}{_{#2}}} }
\newcommand{\PsiHyper}[2]{ {\tensor[_{#1}]{\Psi}{_{#2}}} }
\newcommand{\phiHYPER}[6]{ {\phiHyper{#1}{#2} \!\left(
   \!\!\begin{array}{c} #3 \\ #4 \end{array}\! ;\, #5, \, #6 \! \right)\!} }
\newcommand{\psiHYPER}[6]{ {\psiHyper{#1}{#2} \!\left(
   \!\!\begin{array}{c} #3 \\ #4 \end{array}\! ;\, #5, \, #6 \! \right)} }
\newcommand{\PhiHYPER}[5]{ {\PhiHyper{#1}{#2} \!\left(
   \!\!\begin{array}{c} #3 \\ #4 \end{array}\! ;\, #5 \! \right)\!} }
\newcommand{\PsiHYPER}[5]{ {\PsiHyper{#1}{#2} \!\left(
   \!\!\begin{array}{c} #3 \\ #4 \end{array}\! ;\, #5 \! \right)\!} }
\newcommand{\zerophizero}{ {\phiHyper{0}{0}} }
\newcommand{\ophizero}{ {\phiHyper{1}{0}} }
\newcommand{\zphio}{ {\phiHyper{0}{1}} }
\newcommand{\ophio}{ {\phiHyper{1}{1}} }
\newcommand{\tphio}{ {\phiHyper{2}{1}} }
\newcommand{\tphiz}{ {\phiHyper{2}{0}} }
\newcommand{\tPhio}{ {\PhiHyper{2}{1}} }
\newcommand{\opsio}{ {\psiHyper{1}{1}} }

\def\blue{\textcolor{blue}}
\def\red{\textcolor{red}}
\def\green{\textcolor{green}}
\def\yellow{\textcolor{yellow}}

%
%
\newcommand{\stirlingsubset}[2]{\genfrac{\{}{\}}{0pt}{}{#1}{#2}}
\newcommand{\stirlingcycle}[2]{\genfrac{[}{]}{0pt}{}{#1}{#2}}
\newcommand{\assocstirlingsubset}[3]{{\genfrac{\{}{\}}{0pt}{}{#1}{#2}}_{\! \ge #3}}
\newcommand{\genstirlingsubset}[4]{{\genfrac{\{}{\}}{0pt}{}{#1}{#2}}_{\! #3,#4}}
\newcommand{\irredstirlingsubset}[2]{{\genfrac{\{}{\}}{0pt}{}{#1}{#2}}^{\!\rm irr}}
\newcommand{\euler}[2]{\genfrac{\langle}{\rangle}{0pt}{}{#1}{#2}}
\newcommand{\eulergen}[3]{{\genfrac{\langle}{\rangle}{0pt}{}{#1}{#2}}_{\! #3}}
\newcommand{\eulersecond}[2]{\left\langle\!\! \euler{#1}{#2} \!\!\right\rangle}
\newcommand{\eulersecondgen}[3]{{\left\langle\!\! \euler{#1}{#2} \!\!\right\rangle}_{\! #3}}
\newcommand{\binomvert}[2]{\genfrac{\vert}{\vert}{0pt}{}{#1}{#2}}
\newcommand{\binomsquare}[2]{\genfrac{[}{]}{0pt}{}{#1}{#2}}
\newcommand{\doublebinom}[2]{\left(\!\! \binom{#1}{#2} \!\!\right)}


\newenvironment{sarray}{
             \textfont0=\scriptfont0
             \scriptfont0=\scriptscriptfont0
             \textfont1=\scriptfont1
             \scriptfont1=\scriptscriptfont1
             \textfont2=\scriptfont2
             \scriptfont2=\scriptscriptfont2
             \textfont3=\scriptfont3
             \scriptfont3=\scriptscriptfont3
           \renewcommand{\arraystretch}{0.7}
           \begin{array}{l}}{\end{array}}

\newenvironment{scarray}{
             \textfont0=\scriptfont0
             \scriptfont0=\scriptscriptfont0
             \textfont1=\scriptfont1
             \scriptfont1=\scriptscriptfont1
             \textfont2=\scriptfont2
             \scriptfont2=\scriptscriptfont2
             \textfont3=\scriptfont3
             \scriptfont3=\scriptscriptfont3
           \renewcommand{\arraystretch}{0.7}
           \begin{array}{c}}{\end{array}}


\newcommand*\circled[1]{\tikz[baseline=(char.base)]{
  \node[shape=circle,draw,inner sep=1pt] (char) {#1};}}
\newcommand{\ostar}{{\circledast}}
\newcommand{\ostarN}{{\,\circledast_{\vphantom{\dot{N}}N}\,}}
\newcommand{\ostarPsi}{{\,\circledast_{\vphantom{\dot{\Psi}}\Psi}\,}}
\newcommand{\starN}{{\,\ast_{\vphantom{\dot{N}}N}\,}}
\newcommand{\starpsi}{{\,\ast_{\vphantom{\dot{\bpsi}}\!\bpsi}\,}}
\newcommand{\starone}{{\,\ast_{\vphantom{\dot{1}}1}\,}}
\newcommand{\startwo}{{\,\ast_{\vphantom{\dot{2}}2}\,}}
\newcommand{\starinfty}{{\,\ast_{\vphantom{\dot{\infty}}\infty}\,}}
\newcommand{\starT}{{\,\ast_{\vphantom{\dot{T}}T}\,}}

\newcommand*{\Scale}[2][4]{\scalebox{#1}{$#2$}}

\newcommand*{\Scaletext}[2][4]{\scalebox{#1}{#2}} 

\newcommand{\bolddot}{\boldsymbol{\cdot}}

\clearpage

\tableofcontents

\clearpage

\section{Introduction}

If $(a_n)_{n \ge 0}$ is a sequence of combinatorial numbers or polynomials
with $a_0 = 1$, it is often fruitful to seek to express its
ordinary generating function as a continued fraction of either
Stieltjes type (\textbfit{S-fraction}),
\be
   \sum_{n=0}^\infty a_n t^n
   \;=\;
   \cfrac{1}{1 - \cfrac{\alpha_1 t}{1 - \cfrac{\alpha_2 t}{1 - \cdots}}}
   \label{def.Stype}
   \;\;,
\ee
Thron type (\textbfit{T-fraction}),
\be
   \sum_{n=0}^\infty a_n t^n
   \;=\;
   \cfrac{1}{1 - \delta_1 t - \cfrac{\alpha_1 t}{1 - \delta_2 t - \cfrac{\alpha_2 t}{1 - \cdots}}}
   \label{def.Ttype}
   \;\;,
\ee
or Jacobi type (\textbfit{J-fraction}),
\be
   \sum_{n=0}^\infty a_n t^n
   \;=\;
   \cfrac{1}{1 - \gamma_0 t - \cfrac{\beta_1 t^2}{1 - \gamma_1 t - \cfrac{\beta_2 t^2}{1 - \cdots}}}
   \label{def.Jtype}
   \;\;.
\ee
(Both sides of these expressions are to be interpreted as
formal power series in the indeterminate $t$.)
This line of investigation goes back at least to
Euler \cite{Euler_1760,Euler_1788},
but it gained impetus following Flajolet's \cite{Flajolet_80}
seminal discovery that any S-fraction (resp.~J-fraction)
can be interpreted combinatorially as a generating function
for Dyck (resp.~Motzkin) paths with suitable weights for each rise and fall
(resp.~each rise, fall and level step).
More recently, several authors
\cite{Fusy_15,Oste_15,Josuat-Verges_18,Sokal_totalpos,Elvey-Price-Sokal_wardpoly}
have found a similar combinatorial interpretation
of the general T-fraction:
namely, as a generating function for Schr\"oder paths with suitable
weights for each rise, fall and long level step.
There are now literally dozens of sequences $(a_n)_{n \ge 0}$
of combinatorial numbers or polynomials for which
a continued-fraction expansion of the type
\reff{def.Stype}, \reff{def.Ttype} or \reff{def.Jtype}
is explicitly known.

In a recent paper, Zeng and one of us \cite{Sokal-Zeng_masterpoly}
ran this program in reverse:
starting from a continued fraction in which the coefficients $\balpha$
(or $\bbeta$ and $\bgamma$) contain indeterminates in a nice pattern,
we sought a combinatorial interpretation for the
resulting polynomials $a_n$ --- namely,
as enumerating permutations, set partitions or perfect matchings
according to some natural multivariate statistics.
As a consequence, our results contained many previously obtained identities
as special cases, providing a common refinement of all of them.
In particular, we proved J-fractions enumerating permutations with
10, 18 or infinitely many statistics
that implement the cycle classification of indices (cycle peak, cycle valley,
cycle double rise, cycle double fall, fixed point) together with an
index-refined count of crossings and nestings
(these statistics will be defined in Section~\ref{sec.prelim}).

Subsequently, the two present authors \cite{Deb-Sokal_genocchi}
proved analogous results for D-permu\-ta\-tions
\cite{Lazar_20,Lazar_22,Lazar_23},
which are a subclass of permutations of $[2n]$
(defined in Section~\ref{sec.D-permutations})
that are counted by the Genocchi and median Genocchi numbers:
our T-fractions enumerated D-permutations with 12, 22 or
infinitely many statistics that implement the parity-refined cycle
classification of indices (cycle peak, cycle valley, cycle double rise,
cycle double fall, even fixed point, odd fixed point) together with an
index-refined count of crossings and nestings.
In both papers, we called these results our ``first'' continued fractions.

In both cases, it was natural to try to extend these results
by taking account also of the number of cycles: that is, by including
an additional weight $\lambda^{\cyc(\sigma)}$.  However, it turned out
that it was possible to do so only by renouncing some of the other statistics:
for instance, by counting cycle valleys only with respect to
crossings + nestings, rather than to crossings and nestings separately.
We called these results our ``second'' continued fractions
\cite[Theorems~2.1(b), 2.4, 2.12, 2.14, 2.15]{Sokal-Zeng_masterpoly} 
\cite[Theorems~4.2, 4.7, 4.10]{Deb-Sokal_genocchi}.

Our purpose here is to make a simple but previously overlooked remark:
that in addition to the trivial case $\lambda = 1$, there is one other case
where one need not renounce counting any other statistics, namely,
$\lambda = - 1$.  The reason for this is the following simple lemma,
which relates the number of cycles modulo 2 to the number of
fixed points, cycle peaks (or cycle valleys), and crossings:

\begin{lemma}
   \label{lemma1.1}
Let $\sigma\in \Sym_n$ be a permutation.
Then the following identity holds:
\begin{subeqnarray}
   \cyc & = & \fix + \cpeak + \ucross + \lcross \pmod{2} 
      \slabel{eq.cyc.formula.a} \\[2mm] 
        & = & \fix + \cval + \ucross + \lcross \pmod{2}   \;.
      \slabel{eq.cyc.formula.b}
  \label{eq.cyc.formula}
\end{subeqnarray}
\end{lemma}

\noindent
We will give a precise definition of ucross (number of upper crossings)
and lcross (number of lower crossings) in Section~\ref{subsec.statistics.2},
and then a proof of this lemma in Section~\ref{sec.proofoflemma}.

Using Lemma~\ref{lemma1.1}, it is easy to obtain continued fractions
for the case $\lambda = -1$ as simple corollaries of those for $\lambda = 1$.
That is what we shall do in this paper.

The plan of this paper is as follows:
In Section~\ref{sec.prelim} we give some preliminary definitions
concerning permutation statistics.
In Section~\ref{sec.proofoflemma} we give two proofs of
Lemma~\ref{lemma1.1}:  one topological, and one combinatorial.
Then, in Sections~\ref{sec.permutations} and \ref{sec.D-permutations},
we give our results for permutations and D-permutations, respectively.

Throughout this paper, we shall use two running examples.
The first is the permutation 
\be
\sigma = 9\,3\,7\,4\,6\,11\,2\,8\,10\,1\,5
           = (1,9,10)\,(2,3,7)\,(4)\,(5,6,11)\,(8) \in \Sym_{11};
\ee
the second is the permutation
\begin{eqnarray}
\sigma & = & 7\, 1\, 9\, 2\, 5\, 4\, 8\, 6\, 10\, 3\, 11\, 12\, 14\, 13\,
	\nonumber\\
       & = & (1,7,8,6,4,2)\,(3,9,10)\,(5)\,(11)\,(12)\,(13,14) \in \Sym_{14}.
\end{eqnarray}
We will see later that our second example is a D-permutation.

We remark that since Lemma~\ref{lemma1.1} is a general fact
concerning permutations, it can be applied to {\em any}\/ result
concerning {\em any}\/ subclass of permutations
in which the statistics $\fix$, $\cpeak$ and $\ucross+\lcross$ are handled.

As the reader will have noticed, the present paper builds directly
on the ideas, techniques and intuitions of references
\cite{Sokal-Zeng_masterpoly} and \cite{Deb-Sokal_genocchi}.
Some readers may therefore find it useful to consult those papers first.

\section{Preliminaries}   \label{sec.prelim}


We use the standard notation $[n] \eqdef \{1,\ldots,n\}$.

\subsection{Permutation statistics: The record-and-cycle classification}
     \label{subsec.statistics.1}

Given a permutation $\sigma \in \Sym_N$, an index $i \in [N]$ is called an
\begin{itemize}
   \item {\em excedance}\/ (exc) if $i < \sigma(i)$;
   \item {\em anti-excedance}\/ (aexc) if $i > \sigma(i)$;
   \item {\em fixed point}\/ (fix) if $i = \sigma(i)$.
\end{itemize}
Clearly every index $i$ belongs to exactly one of these three types;
we call this the \textbfit{excedance classification}.
We also say that $i$ is a {\em weak excedance}\/ if $i \le \sigma(i)$,
and a {\em weak anti-excedance}\/ if $i \ge \sigma(i)$.

A more refined classification is as follows:
an index $i \in [N]$ is called a
\begin{itemize}
   \item {\em cycle peak}\/ (cpeak) if $\sigma^{-1}(i) < i > \sigma(i)$;
   \item {\em cycle valley}\/ (cval) if $\sigma^{-1}(i) > i < \sigma(i)$;
   \item {\em cycle double rise}\/ (cdrise) if $\sigma^{-1}(i) < i < \sigma(i)$;
   \item {\em cycle double fall}\/ (cdfall) if $\sigma^{-1}(i) > i > \sigma(i)$;
   \item {\em fixed point}\/ (fix) if $\sigma^{-1}(i) = i = \sigma(i)$.
\end{itemize}
Clearly every index $i$ belongs to exactly one of these five types;
we refer to this classification as the \textbfit{cycle classification}.
Obviously, excedance = cycle valley or cycle double rise,
and anti-excedance = cycle peak or cycle double fall.
We write
\be
   \Cpeak(\sigma)
   \;=\;
   \{ i \colon\: \sigma^{-1}(i) < i > \sigma(i) \}
\ee
for the set of cycle peaks and
\be
   \cpeak(\sigma)  \;=\;  |\Cpeak(\sigma)|
\ee
for its cardinality, and likewise for the others.

On the other hand, an index $i \in [N]$ is called a
\begin{itemize}
   \item {\em record}\/ (rec) (or {\em left-to-right maximum}\/)
         if $\sigma(j) < \sigma(i)$ for all $j < i$
      [note in particular that the indices 1 and $\sigma^{-1}(N)$
       are always records];
   \item {\em antirecord}\/ (arec) (or {\em right-to-left minimum}\/)
         if $\sigma(j) > \sigma(i)$ for all $j > i$
      [note in particular that the indices $N$ and $\sigma^{-1}(1)$
       are always antirecords];
   \item {\em exclusive record}\/ (erec) if it is a record and not also
         an antirecord;
   \item {\em exclusive antirecord}\/ (earec) if it is an antirecord
         and not also a record;
   \item {\em record-antirecord}\/ (rar) (or {\em pivot}\/)
      if it is both a record and an antirecord;
   \item {\em neither-record-antirecord}\/ (nrar) if it is neither a record
      nor an antirecord.
\end{itemize}
Every index $i$ thus belongs to exactly one of the latter four types;
we refer to this classification as the \textbfit{record classification}.
We stress that our records and antirecords are {\em positions}\/, not values.

The record and cycle classifications of indices are related as follows:
\begin{quote}
\begin{itemize}
   \item[(a)]  Every record is a weak excedance,
      and every exclusive record is an excedance.
   \item[(b)]  Every antirecord is a weak anti-excedance,
      and every exclusive antirecord is an anti-excedance.
   \item[(c)]  Every record-antirecord is a fixed point.
\end{itemize}
\end{quote}
Therefore, by applying the record and cycle classifications simultaneously,
we obtain 10 (not~20) disjoint categories \cite{Sokal-Zeng_masterpoly}:
%
%
\medskip
\begin{center}
\begin{tabular}{c|c|c|c|c|c|}
                & cpeak & cval & cdrise & cdfall & fix \\
        \hline
        erec &   & ereccval & ereccdrise & &\\
        earec & eareccpeak & &  & eareccdfall &\\
        rar & & & & & rar \\
        nrar & nrcpeak & nrcval & nrcdrise & nrcdfall & nrfix\\
        \hline
\end{tabular}
\end{center}
\medskip
\noindent
Clearly every index $i$ belongs to exactly one of these 10~types;
we call this the \textbfit{record-and-cycle classification}.

When studying D-permutations, we will use the
\textbfit{parity-refined record-and-cycle classification},
in which we distinguish even and odd fixed points.

\subsubsection{Running example 1}


We consider our first running example
in its cycle notation, \\
$\sigma = (1,9,10)\,(2,3,7)\,(4)\,(5,6,11)\,(8) \in \Sym_{11}.$
The excedance classification of $\sigma$ 
partitions the index set $[11] \eqdef \{1,\ldots,11\}$ as follows:
\be
{\rm Exc} \;=\; \{1, 2, 3, 5, 6, 9\} ,
\qquad {\rm Aexc} \;=\; \{7, 10, 11\} ,
\qquad \Fix \;=\; \{4, 8\} \;.
\label{eq.example.1.excedance.classification}
\ee
Thus, $\exc(\sigma) = 6$, $\aexc(\sigma) = 3$ 
and $\fix(\sigma) = 2$.

Next, we write out the cycle classification of $\sigma$:
\begin{subeqnarray}
	\Cpeak(\sigma) \; = \; \{7, 10, 11 \}  &\qquad&
	\Cval(\sigma) \; = \; \{1,2,5\} \\
	\Cdrise(\sigma) \; = \; \{ 3,6,9\} &\qquad&
        \Cdfall(\sigma) \; = \; \emptyset\\
	\Fix(\sigma) \; = \; \{4, 8\} &&
\label{eq.example.1.cycle.classification}
\end{subeqnarray}
The statistics $\cpeak, \cval, \cdrise, \cdfall$ and $\fix$ are 
simply the cardinalities of these sets, respectively.

For the record classification, we write $\sigma$ as a word,
i.e., $\sigma = 9\:3\:7\:4\:6\:11\:2\:8\:10\:1\:5$.
The~record and antirecord positions are therefore
$\Rec(\sigma) = \{1, 6\}$ and $\Arec(\sigma) = \{10, 11\}$.
The full record classification is
\begin{subeqnarray}
	{\rm Erec}(\sigma) \;=\; \{1,6\} &&
	{\rm Earec}(\sigma) \;=\; \{10,11\} \\
	{\rm Rar}(\sigma) \;=\; \emptyset &&
	{\rm Nrar}(\sigma) \;=\; \{2,3,4,5,7,8,9\}
\label{eq.example.1.record.classification}
\end{subeqnarray}

Finally, the record-and-cycle classification gives us
\begin{subeqnarray}
 	{\rm Eareccpeak}(\sigma) \;=\; \{10, 11\} &&
        {\rm Nrcpeak}(\sigma) \;=\; \{7\} \\
	{\rm Ereccval}(\sigma) \;=\; \{1\} &&
        {\rm Nrcval}(\sigma) \;=\; \{2,5\} \\
	{\rm Erecdrise}(\sigma) \;=\; \{6\} &&
        {\rm Nrcdrise}(\sigma) \;=\; \{3,9\} \\
	{\rm Earecdfall}(\sigma) \;=\; \emptyset &&
        {\rm Nrcdfall}(\sigma) \;=\; \emptyset \\
	{\rm Rar}(\sigma) \;=\; \emptyset &&
	{\rm Nrfix}(\sigma) \; = \; \{4, 8\}
\label{eq.example.1.record.and.cycle.classification}
\end{subeqnarray}

\subsubsection{Running example 2}

We now consider our second running example
in its cycle notation, \\
$\sigma = (1,7,8,6,4,2)\,(3,9,10)\,(5)\,(11)\,(12)\,(13,14) \in \Sym_{14}.$
The excedance classification of $\sigma$ 
partitions the index set $[14] \eqdef \{1,\ldots,14\}$ as follows:
\be
{\rm Exc} \;=\; \{1, 3, 7, 9, 13\} ,
\qquad {\rm Aexc} \;=\; \{2, 4, 6, 8, 10, 14\} ,
\qquad \Fix \;=\; \{5,11,12\} \;.
\label{eq.example.2.excedance.classification}
\ee
Thus, $\exc(\sigma) = 5$, $\aexc(\sigma) = 6$ 
and $\fix(\sigma) = 3$.

Next, we write out the cycle classification of $\sigma$:
\begin{subeqnarray}
	\Cpeak(\sigma) \; = \; \{8, 10, 14 \}  &\qquad&
	\Cval(\sigma) \; = \; \{1, 3, 13\} \\
	\Cdrise(\sigma) \; = \; \{7, 9 \} &\qquad&
	\Cdfall(\sigma) \; = \; \{2, 4, 6\}\\
	\Fix(\sigma) \; = \; \{5,11,12\} &&
\label{eq.example.2.cycle.classification}
\end{subeqnarray}
Once again, the statistics $\cpeak, \cval, \cdrise, \cdfall$ and $\fix$ are 
simply the cardinalities of these sets.

For the record classification, we write $\sigma$ as a word,
$\sigma = 7\: 1\: 9\: 2\: 5\: 4\: 8\: 6\: 10\: 3\: 11\: 12\: 14\: 13$.
The record and antirecord positions are therefore
$\Rec(\sigma) = \{1, 3, 9, 11, 12, 13\}$
and $\Arec(\sigma) = \{2,4, 10, 11, 12, 14\}$.
The full record classification is
\begin{subeqnarray}
	{\rm Erec}(\sigma) \;=\; \{1,3,9,13\} &&
	{\rm Earec}(\sigma) \;=\; \{2,4,10,14\} \\
	{\rm Rar}(\sigma) \;=\; \{11,12\} &&
	{\rm Nrar}(\sigma) \;=\; \{5,6,7,8\}
\label{eq.example.2.record.classification}
\end{subeqnarray}

Finally, the record-and-cycle classification gives us
\begin{subeqnarray}
 	{\rm Eareccpeak}(\sigma) \;=\; \{10,14\} &&
        {\rm Nrcpeak}(\sigma) \;=\; \{8\} \\
	{\rm Ereccval}(\sigma) \;=\; \{1,3,13\} &&
        {\rm Nrcval}(\sigma) \;=\; \emptyset \\
	{\rm Erecdrise}(\sigma) \;=\; \{9\} &&
        {\rm Nrcdrise}(\sigma) \;=\; \{7\} \\
	{\rm Earecdfall}(\sigma) \;=\; \{2,4\} &&
	{\rm Nrcdfall}(\sigma) \;=\; \{6\} \\
	{\rm Rar}(\sigma) \;=\; \{11,12\} &&
	{\rm Nrfix}(\sigma) \; = \; \{5\}
\label{eq.example.2.record.and.cycle.classification}
\end{subeqnarray}

\subsection{Permutation statistics: Crossings and nestings}
     \label{subsec.statistics.2}

We now define (following \cite{Sokal-Zeng_masterpoly})
some permutation statistics that count
\textbfit{crossings} and \textbfit{nestings}.

First we associate to each permutation $\sigma \in \Sym_N$
a pictorial representation
by placing vertices $1,2,\ldots,N$ along a horizontal axis
and then drawing an arc from $i$ to $\sigma(i)$
above (resp.\ below) the horizontal axis
in case $\sigma(i) > i$ [resp.\ $\sigma(i) < i$];
if $\sigma(i) = i$ we do not draw any arc.
Each vertex thus has either
out-degree = in-degree = 1 (if it is not a fixed point) or
out-degree = in-degree = 0 (if it is a fixed point).
Of course, the arrows on the arcs are redundant,
because the arrow on an arc above (resp.\ below) the axis
always points to the right (resp.\ left);
we therefore omit the arrows for simplicity.
See Figures~\ref{fig.pictorial} and \ref{fig.pictorial.2}
for our two running examples.
\begin{figure}[t]
\centering
\vspace*{4cm}
\begin{picture}(60,20)(120, -65)
\setlength{\unitlength}{2mm}
\linethickness{.5mm}
\put(-2,0){\line(1,0){54}}
\put(0,0){\circle*{1,3}}\put(0,0){\makebox(0,-6)[c]{\small 1}}
\put(5,0){\circle*{1,3}}\put(5,0){\makebox(0,-6)[c]{\small 2}}
\put(10,0){\circle*{1,3}}\put(10,0){\makebox(0,-6)[c]{\small 3}}
\put(15,0){\circle*{1,3}}\put(15,0){\makebox(0,-6)[c]{\small 4}}
\put(20,0){\circle*{1,3}}\put(20,0){\makebox(0,-6)[c]{\small 5}}
\put(25,0){\circle*{1,3}}\put(25,0){\makebox(0,-6)[c]{\small 6}}
\put(30,0){\circle*{1,3}}\put(30,0){\makebox(0,-6)[c]{\small 7}}
\put(35,0){ \circle*{1,3}}\put(36,0){\makebox(0,-6)[c]{\small 8}}
\put(40,0){\circle*{1,3}}\put(40,0){\makebox(0,-6)[c]{\small 9}}
\put(45,0){\circle*{1,3}}\put(45,0){\makebox(0,-6)[c]{\small 10}}
\put(50,0){\circle*{1,3}}\put(50,0){\makebox(0,-6)[c]{\small 11}}
\green{\qbezier(0,0)(20,14)(40,0)
\qbezier(40,0)(42.5,6)(45,0)}
\red{\qbezier(4,0)(6.5,5)(9,0)
\qbezier(9,0)(18,10)(29,0)}
\blue{\qbezier(18,0)(20.5,5)(23.5,0)
\qbezier(23.2,0)(36,12)(48.5,0)
\qbezier(18,0)(34,-12)(48.5,0)}
\red{\qbezier(2.5,0)(17,-14)(27.5,0)}
\green{\qbezier(-3,0)(22,-20)(42,0)}
\end{picture}
\caption{
   Pictorial representation of the permutation
   $\sigma = 9\,3\,7\,4\,6\,11\,2\,8\,10\,1\,5
           = (1,9,10)\,(2,3,7)\,(4)\,(5,6,11)\,(8) \in \Sym_{11}$.
 \label{fig.pictorial}
 \vspace*{8mm}
}
\end{figure}

\begin{figure}[t]
\centering
\vspace*{4cm}
\begin{picture}(100,0)(140, -45)
\setlength{\unitlength}{2mm}
\linethickness{.5mm}
\put(-2,0){\line(1,0){69}}
\put(0,0){\circle*{1,3}}\put(0,0){\makebox(0,-6)[c]{\small 1}}
\put(5,0){\circle*{1,3}}\put(5,0){\makebox(0,-6)[c]{\small 2}}
\put(10,0){\circle*{1,3}}\put(10,0){\makebox(0,-6)[c]{\small 3}}
\put(15,0){\circle*{1,3}}\put(15,0){\makebox(0,-6)[c]{\small 4}}
\put(20,0){\circle*{1,3}}\put(20,0){\makebox(0,-6)[c]{\small 5}}
\put(25,0){\circle*{1,3}}\put(25,0){\makebox(0,-6)[c]{\small 6}}
\put(30,0){\circle*{1,3}}\put(30,0){\makebox(0,-6)[c]{\small 7}}
\put(35,0){ \circle*{1,3}}\put(35,0){\makebox(0,-6)[c]{\small 8}}
\put(40,0){\circle*{1,3}}\put(40,0){\makebox(0,-6)[c]{\small 9}}
\put(45,0){\circle*{1,3}}\put(45,0){\makebox(0,-6)[c]{\small 10}}
\put(50,0){\circle*{1,3}}\put(50,0){\makebox(0,-6)[c]{\small 11}}
\put(55,0){\circle*{1,3}}\put(55,0){\makebox(0,-6)[c]{\small 12}}
\put(60,0){\circle*{1,3}}\put(60,0){\makebox(0,-6)[c]{\small 13}}
\put(65,0){\circle*{1,3}}\put(65,0){\makebox(0,-6)[c]{\small 14}}
\green{\qbezier(0,0)(15,16)(30,0)
\qbezier(30,0)(33,6)(36,0)
\qbezier(36,0)(30.5,-8)(25,0)
\qbezier(25,0)(20,-8)(15,0)
\qbezier(15,0)(10,-8)(5,0)
\qbezier(5,0)(2.5,-6)(0,0)
}
\red{\qbezier(9,0)(24,16)(39,0)
\qbezier(39,0)(41.5,6)(44,0)
\qbezier(44,0)(26.5,-18)(9,0)
}
\blue{\qbezier(58.5,0)(61,6)(63.5,0)
\qbezier(63.5,0)(61,-6)(58.5,0)
}
\end{picture}
\caption{
   Pictorial representation of the permutation
   $\sigma = 7\, 1\, 9\, 2\, 5\, 4\, 8\, 6\, 10\, 3\, 11\, 12\, 14\, 13\,
           = (1,7,8,6,4,2)\,(3,9,10)\,(5)\,(11)\,(12)\,(13,14) \in \Sym_{14}$.
   This $\sigma$ is a D-permutation.
 \label{fig.pictorial.2}
 \vspace*{7mm}
}
\end{figure}

We then say that a quadruplet $i < j < k < l$ forms an
\begin{itemize}
   \item \textbfit{upper crossing} (ucross) if $k = \sigma(i)$ and $l = \sigma(j)$;
   \item \textbfit{lower crossing} (lcross) if $i = \sigma(k)$ and $j = \sigma(l)$;
   \item \textbfit{upper nesting}  (unest)  if $l = \sigma(i)$ and $k = \sigma(j)$;
   \item \textbfit{lower nesting}  (lnest)  if $i = \sigma(l)$ and $j = \sigma(k)$.
\end{itemize}
We also consider some ``degenerate'' cases with $j=k$,
by saying that a triplet $i < j < l$ forms an
\begin{itemize}
   \item \textbfit{upper joining} (ujoin) if $j = \sigma(i)$ and $l = \sigma(j)$
      [i.e.\ the index $j$ is a cycle double rise];
   \item \textbfit{lower joining} (ljoin) if $i = \sigma(j)$ and $j = \sigma(l)$
      [i.e.\ the index $j$ is a cycle double fall];
   \item \textbfit{upper pseudo-nesting} (upsnest)
      if $l = \sigma(i)$ and $j = \sigma(j)$;
   \item \textbfit{lower pseudo-nesting} (lpsnest)
      if $i = \sigma(l)$ and $j = \sigma(j)$.
\end{itemize}
These are clearly degenerate cases of crossings and nestings, respectively.
See Figure~\ref{fig.crossnest}.
Note that $\upsnest(\sigma) = \lpsnest(\sigma)$ for all $\sigma$,
since for each fixed point~$j$,
the number of pairs $(i,l)$ with $i < j < l$ such that $l = \sigma(i)$
has to equal the number of such pairs with $i = \sigma(l)$;
we therefore write these two statistics simply as
\be
   \psnest(\sigma) \;\eqdef\; \upsnest(\sigma) \;=\;  \lpsnest(\sigma)
   \;.
\ee
And of course $\ujoin = \cdrise$ and $\ljoin = \cdfall$.

\begin{figure}[p]
\centering
\begin{picture}(30,15)(145, 10)
\setlength{\unitlength}{1.5mm}
\linethickness{.5mm}
\put(2,0){\line(1,0){28}}
\put(5,0){\circle*{1,3}}\put(5,0){\makebox(0,-6)[c]{\small $i$}}
\put(12,0){\circle*{1,3}}\put(12,0){\makebox(0,-6)[c]{\small $j$}}
\put(19,0){\circle*{1,3}}\put(19,0){\makebox(0,-6)[c]{\small $k$}}
\put(26,0){\circle*{1,3}}\put(26,0){\makebox(0,-6)[c]{\small $l$}}
\red{\qbezier(5,0)(12,10)(19,0)}
\blue{\qbezier(11,0)(18,10)(25,0)}
\put(15,-6){\makebox(0,-6)[c]{\small Upper crossing}}
\put(43,0){\line(1,0){28}}
\put(47,0){\circle*{1,3}}\put(47,0){\makebox(0,-6)[c]{\small $i$}}
\put(54,0){\circle*{1,3}}\put(54,0){\makebox(0,-6)[c]{\small $j$}}
\put(61,0){\circle*{1,3}}\put(61,0){\makebox(0,-6)[c]{\small $k$}}
\put(68,0){\circle*{1,3}}\put(68,0){\makebox(0,-6)[c]{\small $l$}}
\red{\qbezier(47,0)(54,-10)(61,0)}
\blue{\qbezier(53,0)(60,-10)(67,0)}
\put(57,-6){\makebox(0,-6)[c]{\small Lower crossing}}
\end{picture}
\\[3.5cm]
\begin{picture}(30,15)(145, 10)
\setlength{\unitlength}{1.5mm}
\linethickness{.5mm}
\put(2,0){\line(1,0){28}}
\put(5,0){\circle*{1,3}}\put(5,0){\makebox(0,-6)[c]{\small $i$}}
\put(12,0){\circle*{1,3}}\put(12,0){\makebox(0,-6)[c]{\small $j$}}
\put(19,0){\circle*{1,3}}\put(19,0){\makebox(0,-6)[c]{\small $k$}}
\put(26,0){\circle*{1,3}}\put(26,0){\makebox(0,-6)[c]{\small $l$}}
\red{\qbezier(5,0)(15.5,10)(26,0)}
\blue{\qbezier(11,0)(14.5,5)(18,0)}
\put(15,-6){\makebox(0,-6)[c]{\small Upper nesting}}
\put(43,0){\line(1,0){28}}
\put(47,0){\circle*{1,3}}\put(47,0){\makebox(0,-6)[c]{\small $i$}}
\put(54,0){\circle*{1,3}}\put(54,0){\makebox(0,-6)[c]{\small $j$}}
\put(61,0){\circle*{1,3}}\put(61,0){\makebox(0,-6)[c]{\small $k$}}
\put(68,0){\circle*{1,3}}\put(68,0){\makebox(0,-6)[c]{\small $l$}}
\red{\qbezier(47,0)(57,-10)(68,0)}
\blue{\qbezier(53,0)(56.5,-5)(60.5,0)}
\put(57,-6){\makebox(0,-6)[c]{\small Lower nesting}}
\end{picture}
\\[3.5cm]
\begin{picture}(30,15)(145, 10)
\setlength{\unitlength}{1.5mm}
\linethickness{.5mm}
\put(2,0){\line(1,0){28}}
\put(5,0){\circle*{1,3}}\put(5,0){\makebox(0,-6)[c]{\small $i$}}
\put(15.5,0){\circle*{1,3}}\put(15.5,0){\makebox(0,-6)[c]{\small $j$}}
\put(26,0){\circle*{1,3}}\put(26,0){\makebox(0,-6)[c]{\small $l$}}
\red{\qbezier(5,0)(10.5,10)(15.5,0)}
\blue{\qbezier(14.75,0)(19.5,10)(25.25,0)}
\put(15,-6){\makebox(0,-6)[c]{\small Upper joining}}
\put(43,0){\line(1,0){28}}
\put(47,0){\circle*{1,3}}\put(47,0){\makebox(0,-6)[c]{\small $i$}}
\put(57.5,0){\circle*{1,3}}\put(57.5,0){\makebox(0,-6)[c]{\small $j$}}
\put(68,0){\circle*{1,3}}\put(68,0){\makebox(0,-6)[c]{\small $l$}}
\red{\qbezier(47,0)(52,-10)(57.25,0)}
\blue{\qbezier(56.75,0)(62,-10)(67,0)}
\put(57,-6){\makebox(0,-6)[c]{\small Lower joining}}
\end{picture}
\\[3.5cm]
\begin{picture}(30,15)(145, 10)
\setlength{\unitlength}{1.5mm}
\linethickness{.5mm}
\put(2,0){\line(1,0){28}}
\put(5,0){\circle*{1,3}}\put(5,0){\makebox(0,-6)[c]{\small $i$}}
\put(15.5,0){\circle*{1,3}}\put(15.5,0){\makebox(0,-6)[c]{\small $j$}}
\put(26,0){\circle*{1,3}}\put(26,0){\makebox(0,-6)[c]{\small $l$}}
\red{\qbezier(5,0)(15.5,10)(26,0)}
\put(15,-6){\makebox(0,-6)[c]{\small Upper pseudo-nesting}}
\put(43,0){\line(1,0){28}}
\put(47,0){\circle*{1,3}}\put(47,0){\makebox(0,-6)[c]{\small $i$}}
\put(57.5,0){\circle*{1,3}}\put(57.5,0){\makebox(0,-6)[c]{\small $j$}}
\put(68,0){\circle*{1,3}}\put(68,0){\makebox(0,-6)[c]{\small $l$}}
\red{\qbezier(47,0)(57.5,-10)(68,0)}
\put(57,-6){\makebox(0,-6)[c]{\small Lower pseudo-nesting}}
\end{picture}
\vspace*{3cm}
\caption{
   Crossing, nesting, joining and pseudo-nesting.
 \label{fig.crossnest}
}
\end{figure}

\begin{figure}[p]
\centering
\begin{picture}(30,15)(145, 10)
\setlength{\unitlength}{1.5mm}
\linethickness{.5mm}
\put(2,0){\line(1,0){28}}
\put(5,0){\circle*{1,3}}\put(5,0){\makebox(0,-6)[c]{\small $i$}}
\put(12,0){\circle*{1,3}}\put(12,0){\makebox(0,-6)[c]{\small $j$}}
\put(19,0){\circle*{1,3}}\put(19,0){\makebox(0,-6)[c]{\small $k$}}
\put(26,0){\circle*{1,3}}\put(26,0){\makebox(0,-6)[c]{\small $l$}}
\red{\qbezier(5,0)(12,10)(19,0)}
\blue{\qbezier(11,0)(18,10)(25,0)}
\blue{\qbezier(10.5,0)(11.75,-1)(13,-2)}
\put(15,-6){\makebox(0,-6)[c]{\small Upper crossing of type cval}}
\put(43,0){\line(1,0){28}}
\put(47,0){\circle*{1,3}}\put(47,0){\makebox(0,-6)[c]{\small $i$}}
\put(54,0){\circle*{1,3}}\put(54,0){\makebox(0,-6)[c]{\small $j$}}
\put(61,0){\circle*{1,3}}\put(61,0){\makebox(0,-6)[c]{\small $k$}}
\put(68,0){\circle*{1,3}}\put(68,0){\makebox(0,-6)[c]{\small $l$}}
\red{\qbezier(47,0)(54,10)(61,0)}
\blue{\qbezier(53,0)(60,10)(67,0)}
\blue{\qbezier(50,2)(51.25,1)(52.5,0)}
\put(57,-6){\makebox(0,-6)[c]{\small Upper crossing of type cdrise}}
\end{picture}
\\[3.5cm]
\begin{picture}(30,15)(145, 10)
\setlength{\unitlength}{1.5mm}
\linethickness{.5mm}
\put(2,0){\line(1,0){28}}
\put(5,0){\circle*{1,3}}\put(5,0){\makebox(0,-6)[c]{\small $i$}}
\put(12,0){\circle*{1,3}}\put(12,0){\makebox(0,-6)[c]{\small $j$}}
\put(19,0){\circle*{1,3}}\put(19,0){\makebox(0,-6)[c]{\small $k$}}
\put(26,0){\circle*{1,3}}\put(26,0){\makebox(0,-6)[c]{\small $l$}}
\red{\qbezier(5,0)(12,-10)(19,0)}
\red{\qbezier(16,2)(17.25, 1)(18.5,0)}
\blue{\qbezier(10.5,0)(18,-10)(25,0)}
\put(15,-6){\makebox(0,-6)[c]{\small Lower crossing of type cpeak}}
\put(43,0){\line(1,0){28}}
\put(47,0){\circle*{1,3}}\put(47,0){\makebox(0,-6)[c]{\small $i$}}
\put(54,0){\circle*{1,3}}\put(54,0){\makebox(0,-6)[c]{\small $j$}}
\put(61,0){\circle*{1,3}}\put(61,0){\makebox(0,-6)[c]{\small $k$}}
\put(68,0){\circle*{1,3}}\put(68,0){\makebox(0,-6)[c]{\small $l$}}
\red{\qbezier(47,0)(54,-10)(61,0)}
\blue{\qbezier(53,0)(60,-10)(67,0)}
\red{\qbezier(59.5,0)(60.75,-1)(62,-2)}
\put(57,-6){\makebox(0,-6)[c]{\small Lower crossing of type cdfall}}
\end{picture}
\\[3.5cm]
\begin{picture}(30,15)(145, 10)
\setlength{\unitlength}{1.5mm}
\linethickness{.5mm}
\put(2,0){\line(1,0){28}}
\put(5,0){\circle*{1,3}}\put(5,0){\makebox(0,-5)[c]{\small $i$}}
\put(12,0){\circle*{1,3}}\put(12,0){\makebox(0,-5)[c]{\small $j$}}
\put(19,0){\circle*{1,3}}\put(19,0){\makebox(0,-5)[c]{\small $k$}}
\put(26,0){\circle*{1,3}}\put(26,0){\makebox(0,-5)[c]{\small $l$}}
\red{\qbezier(5,0)(15.5,10)(26,0)}
\blue{\qbezier(11,0)(14.5,5)(18,0)}
\blue{\qbezier(10.25,0)(11.5,-1)(12.75,-2)}
\put(15,-6){\makebox(0,-6)[c]{\small Upper nesting of type cval}}
\put(43,0){\line(1,0){28}}
\put(47,0){\circle*{1,3}}\put(47,0){\makebox(0,-6)[c]{\small $i$}}
\put(54,0){\circle*{1,3}}\put(54,0){\makebox(0,-6)[c]{\small $j$}}
\put(61,0){\circle*{1,3}}\put(61,0){\makebox(0,-6)[c]{\small $k$}}
\put(68,0){\circle*{1,3}}\put(68,0){\makebox(0,-6)[c]{\small $l$}}
\red{\qbezier(47,0)(57,10)(68,0)}
\blue{\qbezier(53,0)(56.5,5)(60.5,0)}
\blue{\qbezier(50,2)(51.25,1)(52.5,0)}
\put(57,-6){\makebox(0,-6)[c]{\small Upper nesting of type cdrise}}
\end{picture}
\\[3.5cm]
\begin{picture}(30,15)(145, 10)
\setlength{\unitlength}{1.5mm}
\linethickness{.5mm}
\put(2,0){\line(1,0){28}}
\put(5,0){\circle*{1,3}}\put(5,0){\makebox(0,-5)[c]{\small $i$}}
\put(12,0){\circle*{1,3}}\put(12,0){\makebox(0,-5)[c]{\small $j$}}
\put(19,0){\circle*{1,3}}\put(19,0){\makebox(0,-5)[c]{\small $k$}}
\put(26,0){\circle*{1,3}}\put(26,0){\makebox(0,-5)[c]{\small $l$}}
\red{\qbezier(5,0)(15.5,-10)(26,0)}
\blue{\qbezier(11,0)(14.5,-5)(18,0)}
\blue{\qbezier(14.75,2)(16,1)(17.25,0)}
\put(15,-6){\makebox(0,-6)[c]{\small Lower nesting of type cpeak}}
\put(43,0){\line(1,0){28}}
\put(47,0){\circle*{1,3}}\put(47,0){\makebox(0,-5)[c]{\small $i$}}
\put(54,0){\circle*{1,3}}\put(54,0){\makebox(0,-5)[c]{\small $j$}}
\put(61,0){\circle*{1,3}}\put(61,0){\makebox(0,-5)[c]{\small $k$}}
\put(68,0){\circle*{1,3}}\put(68,0){\makebox(0,-5)[c]{\small $l$}}
\red{\qbezier(47,0)(57,-10)(68,0)}
\blue{\qbezier(53,0)(56.5,-5)(60.5,0)}
\blue{\qbezier(59.5,0)(60.75, -1)(62,-2)}
\put(57,-6){\makebox(0,-6)[c]{\small Lower nesting of type cdfall}}
\end{picture}
\vspace*{3cm}
\caption{
   Refined categories of crossing and nesting.
 \label{fig.refined_crossnest}
}
\end{figure}


We can further refine the four crossing/nesting categories
by examining more closely the status of the inner index ($j$ or $k$)
whose {\em outgoing}\/ arc belongs to the crossing or nesting:
that is, $j$ for an upper crossing or nesting,
and $k$ for a lower crossing or nesting:

\medskip
\begin{center}
\begin{tabular}{c|c|c|c|c|}
 & ucross & unest & lcross & lnest \\
\hline
$j\in \Cval$ & ucrosscval & unestcval & &\\
$j\in \Cdrise$ & ucrosscdrise  & unestcdrise & &\\
$k\in \Cpeak$ & & & lcrosscpeak & lnestcpeak\\
$k\in \Cdfall$ & & & lcrosscdfall & lnestcdfall\\
\hline
\end{tabular}
\end{center}
\medskip

\noindent
See Figure~\ref{fig.refined_crossnest}.
Please note that for the ``upper'' quantities
the distinguished index
(i.e.\ the one for which we examine both $\sigma$ and $\sigma^{-1}$)
is in second position ($j$),
while for the ``lower'' quantities
the distinguished index is in third position ($k$).

%

In fact, a central role in our work will be played
(just as in \cite{Sokal-Zeng_masterpoly,Deb-Sokal_genocchi})
by a yet further refinement of these statistics:
rather than counting the {\em total}\/ numbers of quadruplets
$i < j < k < l$ that form upper (resp.~lower) crossings or nestings
of the foregoing types,
we will count the number of upper (resp.~lower) crossings or nestings
that use a particular vertex $j$ (resp.~$k$)
in second (resp.~third) position.
More precisely, we define the
\textbfit{index-refined crossing and nesting statistics}
\begin{subeqnarray}
   \ucross(j,\sigma)
   & = &
   \#\{ i<j<k<l \colon\: k = \sigma(i) \hbox{ and } l = \sigma(j) \}
         \\[2mm]
   \unest(j,\sigma)
   & = &
   \#\{ i<j<k<l \colon\: k = \sigma(j) \hbox{ and } l = \sigma(i) \}
      \\[2mm]
   \lcross(k,\sigma)
   & = &
   \#\{ i<j<k<l \colon\: i = \sigma(k) \hbox{ and } j = \sigma(l) \}
         \\[2mm]
   \lnest(k,\sigma)
   & = &
   \#\{ i<j<k<l \colon\: i = \sigma(l) \hbox{ and } j = \sigma(k) \}
 \label{def.ucrossnestjk}
\end{subeqnarray}
%
%
Note that $\ucross(j,\sigma)$ and $\unest(j,\sigma)$ can be nonzero
only when $j$ is an excedance
(that is, a cycle valley or a cycle double rise),
while $\lcross(k,\sigma)$ and $\lnest(k,\sigma)$ can be nonzero
only when $k$ is an anti-excedance
(that is, a cycle peak or a cycle double fall).

When $j$ is a fixed point, we also define the analogous quantity
for pseudo-nestings:
\be
   \psnest(j,\sigma)
   \;\eqdef\;
   \# \{i < j \colon\:  \sigma(i) > j \}
   \;=\;
   \# \{i > j \colon\:  \sigma(i) < j \}
   \;.
 \label{def.psnestj}
\ee
(Here the two expressions are equal because $\sigma$ is a bijection
 from $[1,j) \cup (j,n]$ to itself.)
In \cite[eq.~(2.20)]{Sokal-Zeng_masterpoly}
this quantity was called the {\em level}\/ of the fixed point $j$
and was denoted $\lev(j,\sigma)$.

\subsubsection{Running example 1}

We first consider our first running example
$\sigma = 9\:3\:7\:4\:6\:11\:2\:8\:10\:1\:5 \\
           = (1,9,10)\,(2,3,7)\,(4)\,(5,6,11)\,(8) \in \Sym_{11}$
and from Figure~\ref{fig.pictorial} 
write out the quadruplets $i<j<k<l$
corresponding to crossings and nestings:
\begin{subeqnarray}
\Ucross(\sigma) & = & \{ 1<6<9<11 ,\; 3<6<7<11  \}\\
        \Lcross(\sigma) & = & \{ 1<5<10<11  ,\;  2<5<7<11 \}\\
        \Unest(\sigma) & = & \{1<2<3<9,\; 1<3<7<9,\; 1<5<6<9,\\
                           && \;\; 3<5<6<7,\; 6<9<10<11\}\\
        \Lnest(\sigma) & = & \{1<2<7<10\}
\label{eq.example.1.cross.nest}
\end{subeqnarray}
We now write out the degenerate cases when $j=k$
but we skip the upper and lower joinings.
The upper and lower pseudo-nestings are:
\begin{subeqnarray}
	\Upsnest(\sigma) & = & \{ 1<4<9,\; 3<4<7,\; 1<8<9  ,\; 6<8<11 \}\\
	\Lpsnest(\sigma) & = & \{ 1<4<10 ,\; 2<4<7 ,\; 1<8<10 ,\; 5<8<11  \}
\label{eq.example.1.upper.lower.pseudonest}
\end{subeqnarray}

Next, we write out the crossings and nestings of $\sigma$
but refined according to the cycle classification 
(which we have already noted down 
in equation~\reff{eq.example.1.cycle.classification}) of
index $j$ for upper crossing or nesting,
and index $k$ for lower crossing or nesting:
\begin{subeqnarray}
	\Ucrosscval(\sigma) & = & \emptyset\\
	\Ucrosscdrise(\sigma) & = & \{ 1<6<9<11 ,\; 3<6<7<11  \}\\
        \Lcrosscpeak(\sigma) & = & \{ 1<5<10<11  ,\;  2<5<7<11 \}\\
        \Lcrosscdfall(\sigma) & = & \emptyset\\
        \Unestcval(\sigma) & = & \{1<2<3<9,\; 1<5<6<9,\; 3<5<6<7\}\\
	\Unestcdrise(\sigma) & = & \{1<3<7<9,\; 6<9<10<11\}\\
        \Lnestcpeak(\sigma) & = & \{1<2<7<10\}\\
        \Lnestcdfall(\sigma) & = & \emptyset
\label{eq.example.1.cross.nest.cycle}
\end{subeqnarray}

Finally, we write out the index-refined crossing and nesting statistics
for $\sigma$.
We make separate tables for the three excedance classes of $\sigma$
[cf.~\reff{eq.example.1.excedance.classification}]:
see Table~\ref{tab.index.refined.1}.

\begin{table}[t]
\begin{center}
\begin{tabular}{|c||*{6}{c|}}
	\hline
	$j\in \Exc(\sigma)$ & 1 & 2 & 3 & 5 & 6 & 9 \\
	\hline
	$\ucross(j,\sigma)$ & 0  & 0  & 0  & 0 & 2 & 0\\
	\hline
	$\unest(j,\sigma)$ & 0  & 1 & 1 & 2 & 0 & 1\\
	\hline
\end{tabular}
\qquad
\begin{tabular}{|c||*{3}{c|}}
	\hline
	$k\in {\rm Aexc}(\sigma)$ & 7 & 10 & 11 \\
	\hline
	$\lcross(k,\sigma)$ & 1  & 1 & 0\\
	\hline
	$\lnest(k,\sigma)$ & 1  & 0 & 0 \\
	\hline
\end{tabular}
	\\[4mm]
\begin{tabular}{|c||*{2}{c|}}
        \hline
        $j\in {\rm Fix}(\sigma)$ & 4 & 8 \\
        \hline
        $\psnest(j,\sigma)$ & 2  & 2 \\
        \hline
\end{tabular}
\end{center}
\vspace*{-3mm}
\caption{Index-refined crossing and nesting statistics for the permutation\\ 
	$\sigma = 9\:3\:7\:4\:6\:11\:2\:8\:10\:1\:5
           = (1,9,10)\,(2,3,7)\,(4)\,(5,6,11)\,(8) \in \Sym_{11}$.}
\vspace*{1cm}
\label{tab.index.refined.1}
\end{table}

\subsubsection{Running example 2}

We now consider our second running example
$\sigma  =  7\: 1\: 9\: 2\: 5\: 4\: 8\: 6\: 10\: 3\: 11\: 12\: 14\: 13\\
        =  (1,7,8,6,4,2)\,(3,9,10)\,(5)\,(11)\,(12)\,(13,14) \in \Sym_{14}$
and from Figure~\ref{fig.pictorial.2} 
write out the quadruplets $i<j<k<l$
corresponding to crossings and nestings:
\begin{subeqnarray}
	 \Ucross(\sigma) & = & \{ 1 < 3 < 7 < 9 \}\\
        \Lcross(\sigma) & = & \{ 2 < 3 < 4 < 10  \}\\
        \Unest(\sigma) & = & \{ 3 < 7 < 8 < 9  \}\\
        \Lnest(\sigma) & = & \{ 3 < 4 < 6 < 10 ,\; 3 < 6 < 8 < 10 \}
\label{eq.example.2.cross.nest}
\end{subeqnarray}
The upper and lower pseudo-nestings are:
\begin{subeqnarray}
	\Upsnest(\sigma) & = & \{  1 < 5 < 7,\; 3< 5< 9\}\\
	\Lpsnest(\sigma) & = & \{ 3< 5< 10,\; 4<5<6  \}
\label{eq.example.2.upper.lower.pseudonest}
\end{subeqnarray}

Next, we write out the crossings and nestings of $\sigma$
but refined according to the cycle classification 
(which we have already noted down 
in equation~\reff{eq.example.2.cycle.classification}) of
index $j$ for upper crossing or nesting,
and index $k$ for lower crossing or nesting:
\begin{subeqnarray}
	\Ucrosscval(\sigma) & = & \{ 1 < 3 < 7 < 9 \}\\
	\Ucrosscdrise(\sigma) & = & \emptyset \\
        \Lcrosscpeak(\sigma) & = & \emptyset\\
        \Lcrosscdfall(\sigma) & = & \{ 2 < 3 < 4 < 10  \}\\
        \Unestcval(\sigma) & = & \emptyset \\
	\Unestcdrise(\sigma) & = & \{ 3 < 7 < 8 < 9  \}\\
        \Lnestcpeak(\sigma) & = & \{ 3 < 6 < 8 < 10 \}\\
        \Lnestcdfall(\sigma) & = & \{ 3 < 4 < 6 < 10 \}
\label{eq.example.2.cross.nest.cycle}
\end{subeqnarray}

Finally, we write out the index-refined crossing and nesting statistics
for $\sigma$.
Again we make separate tables for the three excedance classes of $\sigma$
[cf.~\reff{eq.example.2.excedance.classification}):
see Table~\ref{tab.index.refined.2}.

\begin{table}[t]
\begin{center}
\begin{tabular}{|c||*{5}{c|}}
	\hline
	$j\in \Exc(\sigma)$ & 1 & 3 & 7 & 9 & 13 \\
	\hline
	$\ucross(j,\sigma)$ & 0 & 1 & 0 & 0 & 0\\ 
	\hline
	$\unest(j,\sigma)$ & 0 & 0 & 1 & 0 & 0\\
	\hline
\end{tabular}
\qquad
\begin{tabular}{|c||*{6}{c|}}
	\hline
	$k\in {\rm Aexc}(\sigma)$ & 2 & 4 & 6 & 8 & 10 & 14 \\
	\hline
	$\lcross(k,\sigma)$ & 0 & 1 & 0 & 0 & 0 & 0\\ 
	\hline
	$\lnest(k,\sigma)$ & 0 & 0 & 1 & 1 & 0 & 0\\
	\hline
\end{tabular}
	\\[4mm]
\begin{tabular}{|c||*{3}{c|}}
        \hline
	$j\in {\rm Fix}(\sigma)$ & 5 & 11 & 12\\
        \hline
	$\psnest(j,\sigma)$ & 2 & 0 & 0  \\
        \hline
\end{tabular}
\end{center}
\vspace*{-3mm}
\caption{Index-refined crossing and nesting statistics for the permutation\\
	$\sigma  =  7\: 1\: 9\: 2\: 5\: 4\: 8\: 6\: 10\: 3\: 11\: 12\: 14\: 13
	=  (1,7,8,6,4,2)\,(3,9,10)\,(5)\,(11)\,(12)\,(13,14) \in \Sym_{14}$.}
\vspace*{1cm}
\label{tab.index.refined.2}
\end{table}

\section{Proof of Lemma~\ref{lemma1.1}}  \label{sec.proofoflemma}

We will give two proofs of Lemma~\ref{lemma1.1}:
one topological, and one combinatorial.
The topological proof is extremely satisfying from an intuitive point of view,
but it requires some nontrivial results on the topology of the plane
to make it rigorous.
The combinatorial proof is simple and manifestly rigorous,
but it relies on an identity for the number of inversions
\cite[Lemme~3.1]{DeMedicis_94} \cite[eq.~(40)]{Shin_10}
\cite[Proposition~2.24]{Sokal-Zeng_masterpoly}
whose proof is elementary but not entirely trivial.

\medskip

\begin{figure}[p]
\begin{center}
    \includegraphics{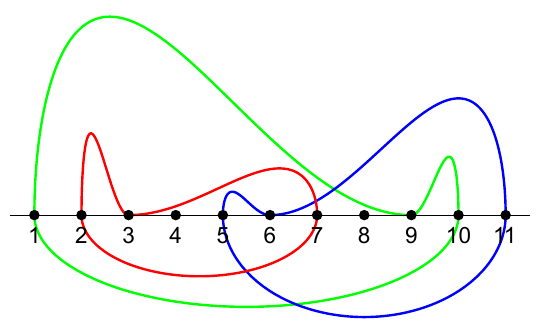}
\end{center}
\vspace*{-6mm}
\caption{Diagram of the permutation \\$\sigma = 9\,3\,7\,4\,6\,11\,2\,8\,10\,1\,5
           = (1,9,10)\,(2,3,7)\,(4)\,(5,6,11)\,(8) \in \Sym_{11}$
   shown in Figure~\ref{fig.pictorial},
   drawn according to the rules stated in the text.
}
\vspace*{10mm}
  \label{fig.pictorial.bis}
\end{figure}

\topologicalproof
Draw the diagram representing the permutation $\sigma$
(Figures~\ref{fig.pictorial} and \ref{fig.pictorial.2})
such that each arc is a $C^1$ non-self-intersecting curve
that has a vertical tangent at each cycle peak and cycle valley
and a horizontal tangent at each cycle double rise and cycle double fall,
and such that each pair of arcs intersects either
zero times (if they do not represent a crossing)
or once transversally (if they do represent a crossing),
and also such that each intersection point involves only two arcs
(see Figures~\ref{fig.pictorial.bis} and \ref{fig.pictorial.2.bis} 
for the examples
 of Figures~\ref{fig.pictorial} and \ref{fig.pictorial.2}, respectively, 
 redrawn according to these rules).
Then each cycle becomes a $C^1$ closed curve
with a finite number of self-intersections,
all of which are transversal double points;
following Whitney \cite[pp.~280--281]{Whitney_37},
we call such a curve {\em normal}\/.
The total number of intersections in the diagram is $\ucross + \lcross$.

\begin{figure}[p]
\centering
\vspace*{4cm}
\begin{picture}(100,0)(140, -45)
\setlength{\unitlength}{2mm}
\linethickness{.5mm}
\green{
\qbezier(1.5,0)(1.5,16)(13,8)
\qbezier(13,8)(26.5,0)(31.5,0)
%
%
\qbezier(31.5,0)(31.7,0)(34.6,0.8)
\qbezier(34.6,0.8)(37.5,1.6)(37,0)
%
%
\qbezier(37,0)(37,-7)(33,-3.5)
\qbezier(33,-3.5)(29,0)(26.5,0)
%
\qbezier(26.5,0)(24,0)(24,-2.5)
\qbezier(24,-2.5)(24,-5)(21.5,-5)
\qbezier(21.5,-5)(19,-5)(19,-2.5)
\qbezier(19,-2.5)(19,0)(16.5,0)
%
\qbezier(16.5,0)(14,0)(14,-2.5)
\qbezier(14,-2.5)(14,-5)(11.5,-5)
\qbezier(11.5,-5)(9,-5)(9,-2.5)
\qbezier(9,-2.5)(9,0)(6.5,0)
%
%
\qbezier(6.5,0)(5.5,0)(3.5,-3.5)
\qbezier(3.5,-3.5)(1.5,-7)(1.5,0)
}
\red{
%
\qbezier(10.5,0)(10.5,20)(23.7,10)
\qbezier(23.7,10)(36.9,0)(40.5,0)
%
%
\qbezier(40.5,0)(41.5,0)(43.5,4)
\qbezier(43.5,4)(45.5,8)(45.5,0)
%
%
\qbezier(45.5,0)(45.5,-8)(28,-8)
\qbezier(28,-8)(10.5,-8)(10.5,0)
}
\blue{\qbezier(60,0)(60,3)(62.5,3)
\qbezier(62.5,3)(65,3)(65,0)
\qbezier(65,0)(65,-3)(62.5,-3)
\qbezier(62.5,-3)(60,-3)(60,0)
}
\put(-2,0){\line(1,0){69}}
\put(0,0){\circle*{1,3}}\put(0,0){\makebox(0,-6)[c]{\small 1}}
\put(5,0){\circle*{1,3}}\put(5,0){\makebox(0,-6)[c]{\small 2}}
\put(10,0){\circle*{1,3}}\put(10,0){\makebox(0,-6)[c]{\small 3}}
\put(15,0){\circle*{1,3}}\put(15,0){\makebox(0,-6)[c]{\small 4}}
\put(20,0){\circle*{1,3}}\put(20,0){\makebox(0,-6)[c]{\small 5}}
\put(25,0){\circle*{1,3}}\put(25,0){\makebox(0,-6)[c]{\small 6}}
\put(30,0){\circle*{1,3}}\put(30,0){\makebox(0,-6)[c]{\small 7}}
\put(35,0){ \circle*{1,3}}\put(35,0){\makebox(0,-6)[c]{\small 8}}
\put(40,0){\circle*{1,3}}\put(40,0){\makebox(0,-6)[c]{\small 9}}
\put(45,0){\circle*{1,3}}\put(45,0){\makebox(0,-6)[c]{\small 10}}
\put(50,0){\circle*{1,3}}\put(50,0){\makebox(0,-6)[c]{\small 11}}
\put(55,0){\circle*{1,3}}\put(55,0){\makebox(0,-6)[c]{\small 12}}
\put(60,0){\circle*{1,3}}\put(60,0){\makebox(0,-6)[c]{\small 13}}
\put(65,0){\circle*{1,3}}\put(65,0){\makebox(0,-6)[c]{\small 14}}

\end{picture}
\caption{
   Diagram of the permutation\\
   $\sigma = 7\, 1\, 9\, 2\, 5\, 4\, 8\, 6\, 10\, 3\, 11\, 12\, 14\, 13\,
           = (1,7,8,6,4,2)\,(3,9,10)\,(5)\,(11)\,(12)\,(13,14) \in \Sym_{14}$
   shown in Figure~\ref{fig.pictorial.2},
   drawn according to the rules stated in the text. 
 \label{fig.pictorial.2.bis}
 \vspace*{7mm}
}
\end{figure}

Each fixed point is of course a cycle.
So we focus henceforth on cycles of length $\ge 2$.
We will prove the following two facts:
\begin{itemize}
   \item[(a)] The number of self-intersections in a cycle
       is equal modulo 2 to the number of cycle peaks
       (or alternatively, cycle valleys) in that cycle, plus 1.
   \item[(b)] The number of intersections between two distinct cycles
       is equal modulo 2 to zero.
\end{itemize}
Together these facts will prove Lemma~\ref{lemma1.1}.

\proofof{(a)}
The {\em rotation angle}\/ (or {\em tangent winding angle}\/)
of a $C^1$ closed curve is the total angle through which
the tangent vector turns while traversing the curve.\footnote{
   In \cite[Section~3]{Umehara_17} and \cite[Chapter~3]{Alencar_22},
   the rotation angle divided by $2\pi$ is called the {\em rotation index}\/.
}
With the above conventions for the arc diagram
(with arcs traversed in the direction of the arrows, i.e.~clockwise)
it is easy to see that the tangent turns by an angle $-\pi$
from each cycle valley to the next cycle peak,
and again by an angle $-\pi$
from each cycle peak to the next cycle valley.
Therefore, a cycle containing $M$ cycle peaks (and hence $M$ cycle valleys)
has a rotation angle $-2\pi M$.
On the other hand, Whitney \cite[Theorem~2]{Whitney_37} proved that
the rotation angle for a $C^1$ normal closed curve $f$ is
\be
   \gamma(f)  \;=\;  2\pi (\mu + N^+ - N^-)
\ee
where $N^+$ (resp.~$N^-$) is the number of positive (resp.~negative) crossings,
and $\mu$ is either $+1$ or $-1$.\footnote{
\begin{samepage}
   The definition of positive and negative crossings
   \cite[p.~281]{Whitney_37} depends on the choice of a starting point
   on the curve;
   if the crossing point is visited first with tangent vector ${\bf v}_1$
   and then with tangent vector ${\bf v}_2$,
   the crossing point is called {\em positive}\/ if
   ${\bf v}_1 \times {\bf v}_2 < 0$ using the right-hand rule,
   and {\em negative}\/ if
   ${\bf v}_1 \times {\bf v}_2 > 0$ using the right-hand rule.
   The hypotheses of \cite[Theorem~2]{Whitney_37}
   require that the starting point be an {\em outside}\/ starting point,
   i.e.~the whole curve must lie on one side
   of the tangent line to the curve at the starting point.
   That requirement is easily fulfilled here,
   e.g.~by taking the starting point to be the
   smallest or largest element of the cycle.
   In this situation, \cite[Theorem~2]{Whitney_37} also specifies
   explicitly whether $\mu$ is $+1$ or $-1$;
   in the present case it is $\mu = -1$.

   See also Umehara and Yamada \cite[pp.~34--38]{Umehara_17}
   for an exposition of Whitney's proof.
   They use the term ``generic'' for what Whitney calls ``normal''.
\end{samepage}
}
It follows that the number of self-intersections in this cycle,
namely $N^+ + N^-$, equals $M+1$ modulo 2.

\proofof{(b)}
This is a general property of $C^1$ normal closed curves in the plane
that have finitely many mutual intersections,
all of which are transversal double points:
in this situation the number of mutual intersections is even.
This intuitively obvious fact goes back at least to
Tait \cite[statement~III]{Tait_1877}.
For completeness we give a proof:

Let $\scrc_1$ and $\scrc_2$ be $C^1$ normal closed curves in the plane;
and suppose that $\scrc_1$ and $\scrc_2$ have finitely many intersections,
all of which are all transversal double points.
Consider first the case in which $\scrc_1$ is a simple closed curve,
i.e.~has no self-intersections.
Then the Jordan Curve Theorem tells us that
$\mathbb{R}^2 \setminus \scrc_1$ has two connected components,
an interior and an exterior.\footnote{
   See e.g.~\cite{Tverberg_80} or \cite[Section~0.3]{Stillwell_93}
   for proofs of the Jordan Curve Theorem.
}
We put an orientation on $\scrc_2$
and traverse $\scrc_2$ from some starting point.
Each time $\scrc_2$ intersects $\scrc_1$,
it must either go from the interior to the exterior of $\scrc_1$
or vice versa (because the intersections are transversal).
Since $\scrc_2$ returns to its starting point,
the number of intersections between $\scrc_2$ and $\scrc_1$ must be even.

When $\scrc_1$ is not a simple closed curve but has finitely many
self-intersections, we can write it as a union of finitely many
simple closed curves $\scrc_1^i$
that are disjoint except for intersections at the self-intersection
points of $\scrc_1$.  (The graph whose vertices are the self-intersection
points and whose edges are the arcs of $\scrc_1$ between two successive
self-intersections is an Eulerian graph;
and an Eulerian graph can be written as the edge-disjoint union
of cycles.)  Then $\scrc_2$ has an even number of intersections with
each $\scrc_1^i$, hence also with $\scrc_1$
(since by hypothesis none of those intersections
occur at the self-intersection points of $\scrc_1$).\footnote{
   Equivalently, the graph $G$ whose vertices are the self-intersection
   points and whose edges are the arcs of $\scrc_1$ between two successive
   self-intersections is Eulerian;
   so its dual $G^*$ is bipartite.
   Then the closed curve $\scrc_2$ must intersect the edges of $G$
   an even number of times.
}
This completes the proof.

\bigskip

This completes the proof of Lemma~\ref{lemma1.1}.
\qed

\combinatorialproof
Let $\cyc(\sigma) = k$,
and let $p_1,\ldots, p_k$ be the sizes of the $k$ cycles of $\sigma$.
Then
\be
   n+k \;=\;  \sum_{i=1}^{k} (p_i+1)
       \;\equiv\;  \#(\text{cycles of $\sigma$ of even length}) \pmod{2} \;.
  \label{eq.nplusk}
\ee
Therefore
\be
   (-1)^{n+k} \;=\; (-1)^{\#(\text{cycles of $\sigma$ of even length})}  \;.
   \label{eq.nplusk.2}
\ee
Here the right-hand side is simply the parity of $\sigma$,
usually denoted $\sgn(\sigma)$.
As is well known (e.g.~\cite[section~7.4]{Loehr_18}),
the parity of $\sigma$ is also given by
\be
   \sgn(\sigma) \;=\; (-1)^{\inv(\sigma)}  \;,
\ee
where
\be
   \inv(\sigma)
   \;\eqdef\;
   \#\{(i,j)\colon\: i < j \textrm{ and } \sigma(i)>\sigma(j) \}
 \label{def.inv}
\ee
is the number of inversions in $\sigma$.
We therefore have
\be
   n+k \;\equiv\; \inv(\sigma) \pmod{2} \;.
 \label{eq.inv.nplusk}
\ee

On the other hand, we recall a formula
\cite[Proposition~2.24]{Sokal-Zeng_masterpoly}
for the number of inversions
in terms of cycle, crossing and nesting statistics:
\be
   \inv \;=\;
   \cval +\cdrise + \cdfall + \ucross +\lcross + 2 (\unest + \lnest + \psnest)
      \;.
\label{eq.invformula}
\ee
Combining \reff{eq.inv.nplusk} and \reff{eq.invformula} yields
\be
   n+k \;\equiv\;  (\cval +\cdrise + \cdfall) + (\ucross +\lcross) \pmod{2}
   \;,
\ee
which can be rewritten as 
\be
   k \;\equiv\; (\cpeak + \fix) + (\ucross +\lcross) \pmod{2}
\ee
since $n =\cpeak+ \cval +\cdrise + \cdfall + \fix$.
This proves \reff{eq.cyc.formula.a}.
Then \reff{eq.cyc.formula.b} follows because $\cpeak = \cval$.
\qed

\subsection{Illustration with examples}

We will verify the various components used in the combinatorial proof 
of Lemma~\ref{lemma1.1} for both of our running examples.

\subsubsection{Running example 1}

First we consider
$\sigma = 9\:3\:7\:4\:6\:11\:2\:8\:10\:1\:5
           = (1,9,10)\,(2,3,7)\,(4)\,(5,6,11)\,(8) \in \Sym_{11}$,
which was depicted in Figure~\ref{fig.pictorial}.
Here $n=11$ and there are $k = 5$ cycles, none of which are of even length.
This confirms~\reff{eq.nplusk} and~\reff{eq.nplusk.2}.

Next we count the number of inversions of $\sigma$.
We record the numbers\\ $\xi_i = {\# \{ j<i \colon\, \sigma(j)>\sigma(i)\}}$,
which are sometimes called the \textbfit{inversion table} of $\sigma$:

\medskip
\begin{center}
\begin{tabular}{|c||*{11}{c|}}
\hline
	$i$ & 1 & 2 & 3 & 4 & 5 & 6 & 7 & 8 & 9 & 10 & 11
	\\
	\hline
	$\sigma(i)$ & 9 & 3 & 7 & 4 & 6 & 11 & 2 & 8 & 10 & 1 & 5
	\\
	\hline
	$\xi_i$ & 0 & 1 & 1 & 2 & 2 & 0 & 6 & 2 & 1  & 9 & 6  \\
	\hline
\end{tabular}
\end{center}
\medskip

%
%

\noindent
Thus we have $\inv(\sigma) = \sum\limits_{i=1}^{11} \xi_i = 30$.
So \reff{eq.inv.nplusk} is also clearly true.

Finally, we will verify equation~\reff{eq.invformula}.
From \reff{eq.example.1.cycle.classification} we obtain the values
\be
\cval(\sigma) \;=\; 3 ,\quad 
\cdrise(\sigma) \;=\; 3  ,\quad
\cdfall(\sigma) \;=\; 0 \;;
\label{eq.cycle.values.example.1}
\ee
and from
\reff{eq.example.1.cross.nest}/\reff{eq.example.1.upper.lower.pseudonest}
we obtain the values
\be
	\ucross(\sigma) \;=\; 2 ,\;
	\lcross(\sigma) \;=\; 2 ,\;
	\unest(\sigma) \;=\; 5 ,\;
	\lnest(\sigma) \;=\; 1 ,\;
	\psnest(\sigma) \;=\; 4 \;.
\label{eq.crossing.values.example.1}
\ee
Using these values we verify \reff{eq.invformula}.

\subsubsection{Running example 2}

Next we consider our second running example
$\sigma = 7\: 1\: 9\: 2\: 5\: 4\: 8\: 6\: 10\: 3\: 11\: 12\: 14\: 13\:
= (1,7,8,6,4,2)\,(3,9,10)\,(5)\,(11)\,(12)\,(13,14) \in \Sym_{14}$,
which was depicted in Figure~\ref{fig.pictorial.2}.
Here $n=14$ and there are $k=6$ cycles, of which two are of even length.
This confirms \reff{eq.nplusk} and~\reff{eq.nplusk.2}.

Next we count the number of inversions of $\sigma$.
We record again the numbers
$\xi_i = \# \{ j<i \:\colon \sigma(j)>\sigma(i)\}$:

\medskip
\begin{center}
\begin{tabular}{|c||*{14}{c|}}
\hline
	$i$ & 1 & 2 & 3 & 4 & 5 & 6 & 7 & 8 & 9 & 10 & 11 & 12 & 13 & 14
        \\
        \hline
	$\sigma(i)$ &
	7 &  1 &  9 &  2 &  5 &  4 &  8 &  6 &  10 &  3 &  11 &  12 &  14 &  13
        \\
        \hline
	$\xi_i$ & 0 & 1  & 0 & 2 & 2 & 3 & 1 & 3 & 0 & 7 & 0 & 0 & 0 & 1  \\
        \hline
\end{tabular}
\end{center}
\medskip

%

\noindent
Thus we have $\inv(\sigma) = \sum\limits_{i=1}^{14} \xi_i = 20$.
So \reff{eq.inv.nplusk} is also clearly true.

Finally, we will verify equation~\reff{eq.invformula}.
From \reff{eq.example.2.cycle.classification} we obtain the values
\be
\cval(\sigma) \;=\; 3 ,\quad
\cdrise(\sigma) \;=\; 2   ,\quad
\cdfall(\sigma) \;=\; 3 \;;
\label{eq.cycle.values.example.2}
\ee
and from
\reff{eq.example.2.cross.nest}/\reff{eq.example.2.upper.lower.pseudonest}
we obtain the values
\be
        \ucross(\sigma) \;=\; 1  ,\;
        \lcross(\sigma) \;=\; 1 ,\;
        \unest(\sigma) \;=\; 1 ,\;
        \lnest(\sigma) \;=\; 2 ,\;
        \psnest(\sigma) \;=\; 2  \;.
\label{eq.crossing.values.example.2}
\ee
Using these values we verify \reff{eq.invformula}.

\section{Results for permutations}   \label{sec.permutations}

We find it convenient to start from the first ``master'' J-fraction
for permutations \cite[Theorem~2.9]{Sokal-Zeng_masterpoly}
and then to specialize.

\subsection{Master J-fraction}

Following \cite[Section~2.7]{Sokal-Zeng_masterpoly},
we introduce five infinite families of indeterminates
$\bsfa = (\sfa_{\ell,\ell'})_{\ell,\ell' \ge 0}$,
$\bsfb = (\sfb_{\ell,\ell'})_{\ell,\ell' \ge 0}$,
$\bsfc = (\sfc_{\ell,\ell'})_{\ell,\ell' \ge 0}$,
$\bsfd = (\sfd_{\ell,\ell'})_{\ell,\ell' \ge 0}$,
$\bsfe = (\sfe_\ell)_{\ell \ge 0}$
and then define the polynomials
\begin{eqnarray}
   & & \hspace*{-6mm}
   P_n(\bsfa,\bsfb,\bsfc,\bsfd,\bsfe,\lambda)
   \;=\;
       \nonumber \\[4mm]
   & &
   \sum_{\sigma \in \Sym_n}
   \lambda^{\cyc(\sigma)}
   \prod\limits_{i \in \Cval(\sigma)}  \! \sfa_{\ucross(i,\sigma),\,\unest(i,\sigma)}
   \prod\limits_{i \in \Cpeak(\sigma)} \!\!  \sfb_{\lcross(i,\sigma),\,\lnest(i,\sigma)}
       \:\times
       \qquad\qquad
       \nonumber \\[1mm]
   & & \qquad\qquad\;\;\,
   \prod\limits_{i \in \Cdfall(\sigma)} \!\!  \sfc_{\lcross(i,\sigma),\,\lnest(i,\sigma)}
   \;
   \prod\limits_{i \in \Cdrise(\sigma)} \!\!  \sfd_{\ucross(i,\sigma),\,\unest(i,\sigma)}
   \, \prod\limits_{i \in \Fix(\sigma)} \sfe_{\psnest(i,\sigma)}
   \;.
   \quad
       \nonumber \\[-3mm]
 \label{def.Qn.firstmaster}
\end{eqnarray}
(This is \cite[eq.~(2.77)]{Sokal-Zeng_masterpoly}
 with a factor $\lambda^{\cyc(\sigma)}$ included.)
Then the first master J-fraction for permutations
\cite[Theorem~2.9]{Sokal-Zeng_masterpoly}
handles the case $\lambda=1$:
it states that the ordinary generating function of the polynomials
$P_n(\bsfa,\bsfb,\bsfc,\bsfd,\bsfe,1)$
has the J-type continued fraction
\begin{eqnarray}
   & & \hspace*{-8mm}
   \sum_{n=0}^\infty P_n(\bsfa,\bsfb,\bsfc,\bsfd,\bsfe,1) \: t^n
   \;=\;
       \nonumber \\
   & & \hspace*{-4mm}
\Scale[0.8]{
   \cfrac{1}{1 - \sfe_0 t - \cfrac{\sfa_{00} \sfb_{00} t^2}{1 -  (\sfc_{00} + \sfd_{00} + \sfe_1) t - \cfrac{(\sfa_{01} + \sfa_{10})(\sfb_{01} + \sfb_{10}) t^2}{1 - (\sfc_{01} + \sfc_{10} + \sfd_{01} + \sfd_{10} + \sfe_2)t - \cfrac{(\sfa_{02} + \sfa_{11} + \sfa_{20})(\sfb_{02} + \sfb_{11} + \sfb_{20}) t^2}{1 - \cdots}}}}
}
       \nonumber \\[1mm]
   \label{eq.thm.permutations.Jtype.final1}
\end{eqnarray}
with coefficients
\begin{subeqnarray}
   \gamma_n  & = &   \biggl( \sum_{\ell=0}^{n-1} \sfc_{\ell,n-1-\ell} \biggr)
               \,+\, \biggl( \sum_{\ell=0}^{n-1} \sfd_{\ell,n-1-\ell} \biggr)
               \,+\, \sfe_n
          \\[1mm]
   \beta_n   & = &   \biggl( \sum_{\ell=0}^{n-1} \sfa_{\ell,n-1-\ell} \biggr)
                  \: \biggl( \sum_{\ell=0}^{n-1} \sfb_{\ell,n-1-\ell} \biggr)
 \label{def.weights.permutations.Jtype.final1}
\end{subeqnarray}

By Lemma~\ref{lemma1.1}, we obtain the case $\lambda=-1$
by inserting a factor $-1$ for each fixed point,
for each cycle peak (or alternatively, cycle valley),
and for each lower or upper crossing.
We therefore have:

\begin{proposition}[Master J-fraction for permutations, $\lambda=-1$]
   \label{prop.permutations.Jtype.final1.lambda=-1}
The ordinary generating function of the polynomials
$P_n(\bsfa,\bsfb,\bsfc,\bsfd,\bsfe,-1)$
has the J-type continued fraction
\begin{eqnarray}
   & & \hspace*{-8mm}
   \sum_{n=0}^\infty P_n(\bsfa,\bsfb,\bsfc,\bsfd,\bsfe,-1) \: t^n
   \;=\;
       \nonumber \\
   & & \hspace*{-4mm}
\Scale[0.8]{
   \cfrac{1}{1 + \sfe_0 t + \cfrac{\sfa_{00} \sfb_{00} t^2}{1 -  (\sfc_{00} + \sfd_{00} - \sfe_1) t + \cfrac{(\sfa_{01} - \sfa_{10})(\sfb_{01} - \sfb_{10}) t^2}{1 - (\sfc_{01} - \sfc_{10} + \sfd_{01} - \sfd_{10} - \sfe_2)t + \cfrac{(\sfa_{02} - \sfa_{11} + \sfa_{20})(\sfb_{02} - \sfb_{11} + \sfb_{20}) t^2}{1 - \cdots}}}}
}
       \nonumber \\[1mm]
\end{eqnarray}
with coefficients
\begin{subeqnarray}
   \gamma_n  & = &   \biggl( \sum_{\ell=0}^{n-1} (-1)^\ell \, \sfc_{\ell,n-1-\ell} \biggr)
               \,+\, \biggl( \sum_{\ell=0}^{n-1} (-1)^\ell \, \sfd_{\ell,n-1-\ell} \biggr)
               \,-\, \sfe_n
          \\[1mm]
   \beta_n   & = & -\, \biggl( \sum_{\ell=0}^{n-1} (-1)^\ell \, \sfa_{\ell,n-1-\ell} \biggr)
                  \: \biggl( \sum_{\ell=0}^{n-1} (-1)^\ell \, \sfb_{\ell,n-1-\ell} \biggr)
 \label{def.weights.permutations.Jtype.final1.lambda=-1}
\end{subeqnarray}
\end{proposition}

We now write out the monomials contributed by our running examples
to the polynomial $P_n(\bsfa,\bsfb,\bsfc,\bsfd,\bsfe,\lambda)$ 
in equation~\reff{def.Qn.firstmaster}
for $n=11$ and $n=14$, respectively.

\subsubsection{Running example 1}

First let us take 
$\sigma = 9\:3\:7\:4\:6\:11\:2\:8\:10\:1\:5
           = (1,9,10)\,(2,3,7)\,(4)\,(5,6,11)\,(8) \in \Sym_{11}$,
which was depicted in Figure~\ref{fig.pictorial}.
Here $n=11$ and $\cyc(\sigma) = 5$.

To obtain the monomial contributed by $\sigma$ in~\reff{def.Qn.firstmaster},
we require the following data for each index $i\in [11]$:
\begin{itemize}
   \item The cycle type of $i$ as per the cycle classification. 
	This determines the letter $\sfa$, $\sfb$, $\sfc$, $\sfd$ or $\sfe$.
	We have already recorded this information
        in~\reff{eq.example.1.cycle.classification}.
   \item The index-refined crossing and nesting statistics for $i$.
        This determines the subscripts $\ell$ and $\ell'$.
        We have already recorded this information
        in Table~\ref{tab.index.refined.1}.
\end{itemize}
We copy these data into the following table:

\medskip
\begin{center}
\begin{tabular}{|c||*{11}{c|}}
\cline{2-12}
\multicolumn{1}{c||}{\quad}
         & 1 & 2 & 3 & 4 & 5 & 6 & 7 & 8 & 9 & 10 & 11
        \\
        \hline
        Letter & \sfa  & \sfa & \sfd  & \sfe  & \sfa  & \sfd  & \sfb  & \sfe  & \sfd  & \sfb & \sfb 
        \\
        \hline
        First subscript & 0 & 0 & 0 & 2 & 0 & 2 & 1 & 2 & 0  & 1 & 0  \\
        \hline
        Second subscript & 0 & 1 & 1 &  & 2 & 0 & 1 &  & 1  & 0 & 0  \\
        \hline
\end{tabular}
\end{center}
\medskip

\noindent
We therefore see that the monomial contributed
to $P_n(\bsfa,\bsfb,\bsfc,\bsfd,\bsfe,\lambda)$
by this particular permutation $\sigma$ is
\be
\lambda^5\,
\sfa_{0,0} \, \sfa_{0,1} \, \sfa_{0,2} \,
\sfb_{0,0} \, \sfb_{1,0} \, \sfb_{1,1} \,
\sfd_{0,1}^2 \, \sfd_{2,0} \,
\sfe_{2}^2\;.
\ee

\subsubsection{Running example 2}

We now consider our second running example
$\sigma  =  7\: 1\: 9\: 2\: 5\: 4\: 8\: 6\: 10\: 3\: 11\: 12\: 14\: 13\\
        =  (1,7,8,6,4,2)\,(3,9,10)\,(5)\,(11)\,(12)\,(13,14) \in \Sym_{14}$,
which was depicted in Figure~\ref{fig.pictorial.2}.
Here $n=14$ and $\cyc(\sigma) = 6$.

To obtain the monomial contributed by $\sigma$ in~\reff{def.Qn.firstmaster},
we again copy the required data
from equation~\reff{eq.example.2.cycle.classification}
and Table~\ref{tab.index.refined.2}:

\medskip
\begin{center}
\begin{tabular}{|c||*{14}{c|}}
\cline{2-15}
\multicolumn{1}{c||}{\quad}
	 & 1 & 2 & 3 & 4 & 5 & 6 & 7 & 8 & 9 & 10 & 11 & 12 & 13 & 14
        \\
        \hline
	Letter & \sfa  & \sfc  & \sfa  & \sfc  & \sfe  & \sfc  & \sfd  & \sfb  & \sfd  & \sfb  & \sfe & \sfe & \sfa & \sfb
        \\
        \hline
	First subscript & 0 & 0 & 1 & 1 & 2 & 0 & 0 & 0 & 0  & 0 & 0 & 0 & 0 & 0  \\
        \hline
	Second subscript & 0 & 0 & 0 & 0  &  & 1 & 1 & 1 & 0 & 0 &  & & 0 & 0 \\
        \hline
\end{tabular}
\end{center}
\medskip

\noindent
We therefore see that the monomial contributed
to $P_n(\bsfa,\bsfb,\bsfc,\bsfd,\bsfe,\lambda)$
by this particular permutation $\sigma$ is
\be
\lambda^6\,
\sfa_{0,0}^2 \, \sfa_{1,0} \,
\sfb_{0,0}^2 \, \sfb_{0,1} \,
\sfc_{0,0} \, \sfc_{0,1} \, \sfc_{1,0} \,
\sfd_{0,0} \, \sfd_{0,1} \,
\sfe_{0}^2 \, \sfe_{2}\;.
\label{eq.monomial.example.2}
\ee

\subsection[\texorpdfstring{$p,q$}{p,q} J-fraction]{{\texorpdfstring{$\bm{p,q}$}{p,q}} J-fraction}

Consider now the polynomial \cite[eq.~(2.92)]{Sokal-Zeng_masterpoly}
\begin{eqnarray}
   & &
   P_n(x_1,x_2,y_1,y_2,u_1,u_2,v_1,v_2,\bw,p_{+1},p_{+2},p_{-1},p_{-2},q_{+1},q_{+2},q_{-1},q_{-2},s,\lambda)
   \;=\;
       \nonumber \\[4mm]
   & & \qquad
   \sum_{\sigma \in \Sym_n}
   x_1^{\eareccpeak(\sigma)} x_2^{\eareccdfall(\sigma)} 
   y_1^{\ereccval(\sigma)} y_2^{\ereccdrise(\sigma)} 
   \:\times
       \qquad\qquad
       \nonumber \\[-1mm]
   & & \qquad\qquad\:
   u_1^{\nrcpeak(\sigma)} u_2^{\nrcdfall(\sigma)}
   v_1^{\nrcval(\sigma)} v_2^{\nrcdrise(\sigma)}
   \bw^{\bfix(\sigma)}  \:\times
       \qquad\qquad
       \nonumber \\[3mm]
   & & \qquad\qquad\:
   p_{+1}^{\ucrosscval(\sigma)}
   p_{+2}^{\ucrosscdrise(\sigma)}
   p_{-1}^{\lcrosscpeak(\sigma)}
   p_{-2}^{\lcrosscdfall(\sigma)}
          \:\times
       \qquad\qquad
       \nonumber \\[3mm]
   & & \qquad\qquad\:
   q_{+1}^{\unestcval(\sigma)}
   q_{+2}^{\unestcdrise(\sigma)}
   q_{-1}^{\lnestcpeak(\sigma)}
   q_{-2}^{\lnestcdfall(\sigma)}
   s^{\psnest(\sigma)}
   \lambda^{\cyc(\sigma)}
   \;,
 \label{def.Qn.BIG.cyc.pq}
\end{eqnarray}
where the various statistics have been defined in
\cite[Sections~2.3 and 2.5]{Sokal-Zeng_masterpoly}.
In order to distinguish records and antirecords,
we use the following general fact about permutations
\cite[Lemma~2.10]{Sokal-Zeng_masterpoly}:
\begin{itemize}
   \item[(a)]  If $i$ is a cycle valley or cycle double rise,
       then $i$ is a record if and only if $\unest(i,\sigma) = 0$;
       and in this case it is an exclusive record.
   \item[(b)]  If $i$ is a cycle peak or cycle double fall,
       then $i$ is an antirecord if and only if $\lnest(i,\sigma) = 0$;
       and in this case it is an exclusive antirecord.
\end{itemize}
It follows that the polynomial \reff{def.Qn.BIG.cyc.pq}
is obtained from \reff{def.Qn.firstmaster}
by making the specializations \cite[eq.~(2.81)]{Sokal-Zeng_masterpoly}
\begin{subeqnarray}
   \sfa_{\ell,\ell'}
   & = &
   p_{+1}^\ell q_{+1}^{\ell'} \,\times\,
   \begin{cases}
       y_1    &  \textrm{if $\ell' = 0$}  \\
       v_1    &  \textrm{if $\ell' \ge 1$}
   \end{cases}
       \\[2mm]
   \sfb_{\ell,\ell'}
   & = &
   p_{-1}^\ell q_{-1}^{\ell'} \,\times\,
   \begin{cases}
       x_1    &  \textrm{if $\ell' = 0$}  \\
       u_1    &  \textrm{if $\ell' \ge 1$}
   \end{cases}
       \\[2mm]
   \sfc_{\ell,\ell'}
   & = &
   p_{-2}^\ell q_{-2}^{\ell'} \,\times\,
   \begin{cases}
       x_2    &  \textrm{if $\ell' = 0$}  \\
       u_2    &  \textrm{if $\ell' \ge 1$}
   \end{cases}
       \\[2mm]
   \sfd_{\ell,\ell'}
   & = &
   p_{+2}^\ell q_{+2}^{\ell'} \,\times\,
   \begin{cases}
       y_2    &  \textrm{if $\ell' = 0$}  \\
       v_2    &  \textrm{if $\ell' \ge 1$}
   \end{cases}
       \\[2mm]
   \sfe_\ell  & = &  s^\ell w_\ell
 \label{eq.Qn.BIG.specializations}
\end{subeqnarray}
Making these specializations in
Proposition~\ref{prop.permutations.Jtype.final1.lambda=-1}
--- or equivalently, attaching a minus sign to the variables
$x_1, u_1, p_{+1}, p_{+2}, p_{-1}, p_{-2}, w_i$
in \cite[Theorem~2.7]{Sokal-Zeng_masterpoly} ---
we obtain:

\begin{proposition}[$p,q$ J-fraction for permutations, $\lambda=-1$]
   \label{prop.permutations.Jtype.pq.lambda=-1}
The ordinary generating function of the polynomials
\reff{def.Qn.BIG.cyc.pq} at $\lambda = -1$
has the J-type continued fraction
\begin{eqnarray}
   & & \hspace*{-7mm}
   \sum_{n=0}^\infty
       P_n(x_1,x_2,y_1,y_2,u_1,u_2,v_1,v_2,\bw,p_{+1},p_{+2},p_{-1},p_{-2},q_{+1},q_{+2},q_{-1},q_{-2},s,-1) \: t^n
   \;=\;
       \nonumber \\
   & & \!\!\!\!\!\!
\Scale[0.65]{
   \cfrac{1}{1 + w_0 t + \cfrac{x_1 y_1 t^2}{1 -  (x_2\!+\!y_2\!-\!sw_1) t + \cfrac{(-p_{-1} x_1\!+\! q_{-1} u_1)(-p_{+1} y_1\!+\! q_{+1} v_1) t^2}{1 - (-p_{-2} x_2\!+\!q_{-2} u_2\!-\!p_{+2} y_2\!+\! q_{+2} v_2\!-\! s^2 w_2)t + \cfrac{(p_{-1}^2 x_1\!+\! q_{-1} [2]_{-p_{-1},q_{-1}} u_1)(p_{+1}^2 y_1\!+\! q_{+1} [2]_{-p_{+1},q_{+1}} v_1) t^2}{1 - \cdots}}}}
}
       \nonumber \\[1mm]
   \label{eq.thm.perm.pq.Jtype.BIG}
\end{eqnarray}
with coefficients
\begin{subeqnarray}
   \gamma_0  & = &   -w_0
      \slabel{def.weights.perm.pq.Jtype.BIG.a}  \\[1mm]
   \gamma_n  & = &
\Scale[0.87]{
    ((-p_{-2})^{n-1} x_2 + q_{-2} \, [n-1]_{-p_{-2},q_{-2}} u_2)
\:+\: ((-p_{+2})^{n-1} y_2 + q_{+2} \, [n-1]_{-p_{+2},q_{+2}} v_2) \:-\: s^n w_n
}
       \nonumber  \\
    & & \hspace*{3.7in}
        \quad\hbox{for $n \ge 1$}
      \slabel{def.weights.perm.pq.Jtype.BIG.b}  \\
   \beta_n   & = &
\Scale[0.87]{
    - \, ((-p_{-1})^{n-1} x_1 + q_{-1} \, [n-1]_{-p_{-1},q_{-1}} u_1)
           \: ((-p_{+1})^{n-1} y_1 + q_{+1} \, [n-1]_{-p_{+1},q_{+1}} v_1)
}
      \slabel{def.weights.perm.pq.Jtype.BIG.c}
 \label{def.weights.perm.pq.Jtype.BIG}
\end{subeqnarray}
\end{proposition}

\begin{sloppypar}
We now write out the monomials contributed by our running examples
to the polynomial 
$P_n(x_1,x_2,y_1,y_2,u_1,u_2,v_1,v_2,\bw,p_{+1},p_{+2},p_{-1},p_{-2},q_{+1},q_{+2},q_{-1},q_{-2},s,\lambda)$
in equation~\reff{def.Qn.BIG.cyc.pq}
for $n=11$ and $n=14$, respectively.
\end{sloppypar}

\subsubsection{Running example 1}

First let us take 
$\sigma = 9\:3\:7\:4\:6\:11\:2\:8\:10\:1\:5
           = (1,9,10)\,(2,3,7)\,(4)\,(5,6,11)\,(8) \in \Sym_{11}$,
which was depicted in Figure~\ref{fig.pictorial}.
Here $n=11$ and $\cyc(\sigma) = 5$.

To obtain the monomial contributed to \reff{def.Qn.BIG.cyc.pq} by $\sigma$,
we require the following data for each index $i\in [11]$:
\begin{itemize}
   \item The cycle-and-record  type of $i$ as per the cycle-and-record
		classification. 
		This determines the letter $x$, $y$, $u$ or $v$
              along with the subscript $1$ or $2$. 
		We have already recorded this information
		in~\reff{eq.example.1.record.and.cycle.classification}.
\end{itemize}
We also require the total numbers
of crossings and nestings refined according to cycle type. 
We have already recorded this information in
\reff{eq.example.1.upper.lower.pseudonest}/%
\reff{eq.example.1.cross.nest.cycle}.
Copying all these data together, we find that
the monomial contributed to \reff{def.Qn.BIG.cyc.pq}
by the permutation $\sigma$ is
\be
\lambda^5\,
x_1^2\,
y_1 \,
y_2 \,
u_1 \,
v_1^2 \,
v_2^2 \,
w_2^2 \,\,
p_{+2}^2 \,
p_{-1}^2 \,
q_{+1}^3 \,
q_{+2}^2 \,
q_{-1} \,
s^4
   \;.
\ee

\subsubsection{Running example 2}

We now consider our second running example
$\sigma  =  7\: 1\: 9\: 2\: 5\: 4\: 8\: 6\: 10\: 3\: 11\: 12\: 14\: 13\\
        =  (1,7,8,6,4,2)\,(3,9,10)\,(5)\,(11)\,(12)\,(13,14) \in \Sym_{14}$,
which was depicted in Figure~\ref{fig.pictorial.2}.
Here $n=14$ and $\cyc(\sigma) = 6$.

Copying the required data from
\reff{eq.example.2.record.and.cycle.classification}/%
\reff{eq.example.2.upper.lower.pseudonest}/%
\reff{eq.example.2.cross.nest.cycle},
we find that the monomial contributed to \reff{def.Qn.BIG.cyc.pq}
by the permutation $\sigma$ is
\be
\lambda^6 \:
x_1^2  \:
x_2^2 \:
y_1^3 \:
y_2 \:
u_1 \:
u_2 \:
v_2 \:
w_0^2 \:
w_2 \:
p_{+1} \: 
p_{-2} \:
q_{+2} \:
q_{-1} \:
q_{-2} \:
s^2  \;.
\label{eq.monomial.example.2.pq}
\ee

\subsection{Simple J-fraction}

And finally, we can obtain the polynomials without crossing and nesting
statistics,
\begin{eqnarray}
   & &
   \hspace*{-6mm}
   P_n(x_1,x_2,y_1,y_2,u_1,u_2,v_1,v_2,\bw,\lambda)
   \;=\;
       \nonumber \\[4mm]
   & & \qquad
   \sum_{\sigma \in \Sym_n}
   x_1^{\eareccpeak(\sigma)} x_2^{\eareccdfall(\sigma)} 
   y_1^{\ereccval(\sigma)} y_2^{\ereccdrise(\sigma)} 
   \:\times
       \qquad\qquad
       \nonumber \\[-1mm]
   & & \qquad\qquad\:
   u_1^{\nrcpeak(\sigma)} u_2^{\nrcdfall(\sigma)}
   v_1^{\nrcval(\sigma)} v_2^{\nrcdrise(\sigma)}
   \bw^{\bfix(\sigma)}
   \lambda^{\cyc(\sigma)}
   \;,
 \label{def.Qn.BIG.cyc}
\end{eqnarray}
by setting
$p_{+1} = p_{+2} = p_{-1} = p_{-2} = q_{+1} = q_{+2} = q_{-1} = q_{-2} = s = 1$
in \reff{def.Qn.BIG.cyc.pq}.
Making this same specialization in
Proposition~\ref{prop.permutations.Jtype.pq.lambda=-1}
and observing that
\be
   [n-1]_{-1,1}  \;=\;
   \begin{cases}
       1 & \textrm{if $n$ is even} \\
       0 & \textrm{if $n$ is odd}
   \end{cases}
   \label{eq.n.minus.one.pq}
\ee
we obtain:

\begin{proposition}[Simple J-fraction for permutations, $\lambda=-1$]
   \label{prop.permutations.Jtype.lambda=-1}
The ordinary generating function of the polynomials
\reff{def.Qn.BIG.cyc} at $\lambda = -1$
has the J-type continued fraction
\begin{eqnarray}
   & & \hspace*{-7mm}
   \sum_{n=0}^\infty
       P_n(x_1,x_2,y_1,y_2,u_1,u_2,v_1,v_2,\bw,-1) \: t^n
   \;=\;
       \nonumber \\
   & & \!\!\!\!
   \cfrac{1}{1 + w_0 t + \cfrac{x_1 y_1 t^2}{1 -  (x_2\!+\!y_2\!-\!w_1) t + \cfrac{(x_1\!-\! u_1)(y_1\!-\! v_1) t^2}{1 - (-x_2\!+\! u_2\!-\! y_2\!+\! v_2\!-\! w_2)t + \cfrac{x_1 y_1 t^2}{1 - \cdots}}}}
       \nonumber \\[1mm]
   \label{eq.thm.perm.Jtype.BIG}
\end{eqnarray}
with coefficients
\begin{subeqnarray}
   \gamma_0  & = &   - w_0
      \slabel{def.weights.perm.Jtype.BIG.a}  \\[1mm]
   \gamma_n  & = &
      \begin{cases}
        x_2 + y_2 - w_n         & \textrm{if $n$ is odd}  \\
        -x_2 + u_2 - y_2 + v_2 - w_n   & \textrm{if $n$ is even and $\ge 2$}
      \end{cases}
      \slabel{def.weights.perm.Jtype.BIG.b}  \\[1mm]
   \beta_n   & = &
      \begin{cases}
        - x_1 y_1                    & \textrm{if $n$ is odd}  \\
        - (x_1 - u_1) (y_1 - v_1)    & \textrm{if $n$ is even}
      \end{cases}
      \slabel{def.weights.perm.Jtype.BIG.c}
 \label{def.weights.perm.Jtype.BIG}
\end{subeqnarray}
\end{proposition}

\subsection{Corollary for cycle-alternating permutations}

We recall \cite{Deutsch_11, Sokal-Zeng_masterpoly, Deb-Sokal_cyclealt}
that a \textbfit{cycle-alternating permutation} is a permutation of $[2n]$
that has no cycle double rises, cycle double falls, or fixed points;
Deutsch and Elizalde \cite[Proposition~2.2]{Deutsch_11}
showed that the number of cycle-alternating permutations of $[2n]$
is the secant number $E_{2n}$
(see also Dumont \cite[pp.~37, 40]{Dumont_86}
 and Biane \cite[section~6]{Biane_93}).
In this subsection, 
we will obtain continued fractions for cycle-alternating permutations
at $\lambda=-1$ by specializing our master J-fraction
(Proposition~\ref{prop.permutations.Jtype.final1.lambda=-1})
to suppress cycle double rises, cycle double falls and fixed points,
and then using \cite[Lemma~4.2]{Deb-Sokal_cyclealt}
to interpret the parity of cycle peaks and cycle valleys
in terms of crossings and nestings.

Let $P_n(\bsfa,\bsfb, \lambda)$ denote the polynomial
\reff{def.Qn.firstmaster}
specialized to $\bsfc = \bsfd = \bsfe = \bzero$;
it enumerates cycle-alternating permutations
according to the index-refined crossing and nesting statistics
associated to its cycle peaks and cycle valleys.
Note that $P_n$ is nonvanishing only for even $n$.
The J-fraction of Proposition~\ref{prop.permutations.Jtype.final1.lambda=-1}
then becomes an S-fraction in the variable $t^2$;
after changing $t^2$ to $t$, we have:

\begin{proposition}[Master S-fraction for cycle-alternating permutations, $\lambda=-1$]
   \label{prop.cyclealt.Stype.final1.lambda=-1}
The ordinary generating function of the polynomials
$P_{2n}(\bsfa,\bsfb,-1)$
has the S-type continued fraction
	\begin{eqnarray}
   \sum_{n=0}^\infty P_{2n}(\bsfa,\bsfb,-1) \: t^n
   \;=\;
   \cfrac{1}{1 + \cfrac{\sfa_{00} \sfb_{00} t}{1 + \cfrac{(\sfa_{01} - \sfa_{10})(\sfb_{01} - \sfb_{10}) t}{1  + \cfrac{(\sfa_{02} - \sfa_{11} + \sfa_{20})(\sfb_{02} - \sfb_{11} + \sfb_{20}) t}{1 - \cdots}}}} \nonumber \\
\end{eqnarray}
with coefficients
\be
\alpha_n   \; = \;  -\, \biggl( \sum_{\ell=0}^{n-1} (-1)^\ell \, \sfa_{\ell,n-1-\ell} \biggr) \: \biggl( \sum_{\ell=0}^{n-1} (-1)^\ell \, \sfb_{\ell,n-1-\ell} \biggr)
   \;.
 \label{def.weights.cyclalt.Stype.lambda=-1}
\ee
\end{proposition}

We can use this master S-fraction to obtain a continued fraction
that distinguishes cycle peaks and cycle valleys according to their parity.
To do this, we use \cite[Lemma~4.2]{Deb-Sokal_cyclealt}:

\begin{lemma}[Key lemma from \protect\cite{Deb-Sokal_cyclealt}]
   \label{lemma.cycle-alt}
If $\sigma$ is a cycle-alternating permutation of $[2n]$, then
\begin{subeqnarray}
   \hbox{\rm cycle valleys:}
   & \!\!  &
   \ucross(i,\sigma) \,+\, \unest(i,\sigma)
   \;=\;
   i-1 \: \pmod{2}
         \\[1mm]
   \hbox{\rm cycle peaks:}
   & \!\!  &
   \lcross(i,\sigma) \,+\, \lnest(i,\sigma)
   \;=\;
   i \: \pmod{2}
\end{subeqnarray}
for all $i \in [2n]$.
\end{lemma}

Consider now the polynomials
\begin{eqnarray}
   & & \hspace*{-6mm}
	Q_n(\xe,\ye,\ue,\ve,\xo,\yo,\uo,\vo,p_{-1},p_{-2},p_{+1},p_{+2},q_{-1},q_{-2},q_{+1},q_{+2},\lambda)
   \;=\;
        \nonumber \\[4mm]
   & & \qquad
   \sum_{\sigma \in \Sym^{\rm ca}_{2n}}
   \xe^{\eareccpeakeven(\sigma)}
   \ye^{\ereccvaleven(\sigma)}
   \ue^{\nrcpeakeven(\sigma)}
   \ve^{\nrcvaleven(\sigma)}
          \:\times
       \qquad\qquad
       \nonumber \\[-2mm]
   & & \qquad\qquad\;\,
   \xo^{\eareccpeakodd(\sigma)}
   \yo^{\ereccvalodd(\sigma)}
   \uo^{\nrcpeakodd(\sigma)}
   \vo^{\nrcvalodd(\sigma)}
	\:\times
       \qquad\qquad
       \nonumber \\[2mm]
   & & \qquad\qquad\;\,
   p_{-1}^{{\rm lcrosscpeakeven}(\sigma)}
   p_{-2}^{{\rm lcrosscpeakodd}(\sigma)}
   p_{+1}^{{\rm ucrosscvalodd}(\sigma)}
   p_{+2}^{{\rm ucrosscvaleven}(\sigma)}
	 \:\times
       \qquad\qquad
       \nonumber \\[2mm]
   & & \qquad\qquad\;\,
   q_{-1}^{{\rm lnestcpeakeven}(\sigma)}
   q_{-2}^{{\rm lnestcpeakodd}(\sigma)}
   q_{+1}^{{\rm unestcvalodd}(\sigma)}
   q_{+2}^{{\rm unestcvaleven}(\sigma)}
   \lambda^{\cyc(\sigma)}
   \;,
\label{def.Qn.cycle-alternating.evenodd.pqgen}
\end{eqnarray}
where the various statistics are defined as
\be
   \eareccpeakeven(\sigma)
   \;=\;
   |\Eareccpeakeven(\sigma)|
   \;=\;
   |\Arec(\sigma) \cap \Cpeak(\sigma) \cap \Even|
\ee
\be
{\rm lcrosscpeakeven}(\sigma) \;=\; \sum_{k\in \Cpeak(\sigma) \,\cap \,{\rm Even}} \lcross(k,\sigma)
\ee
and likewise for the others.
These polynomials are the same as the polynomials 
\cite[eq.~(4.29)]{Deb-Sokal_cyclealt}
except for the extra factor $\lambda^{\cyc(\sigma)}$.
As before, we use \cite[Lemma~2.10]{Sokal-Zeng_masterpoly}
to distinguish records and antirecords;
we also use Lemma~\ref{lemma.cycle-alt}
to distinguish cycle peaks and cycle valleys according to their parity.
It follows that
the polynomials \reff{def.Qn.cycle-alternating.evenodd.pqgen}
can be obtained from the master polynomials $P_n(\bsfa,\bsfb, \lambda)$
by making the specializations \cite[eq.~(4.33)]{Deb-Sokal_cyclealt}
\begin{subeqnarray}
   \sfa_{\ell,\ell'}
   & = &
   \begin{cases}
       p_{+1}^{\ell} \yo    &  \textrm{if $\ell' = 0$ and $\ell+\ell'$ is even}  \\
       p_{+1}^{\ell} q_{+1}^{\ell'}\vo    &  \textrm{if $\ell' \ge 1$ and $\ell+\ell'$ is even} \\
       p_{+2}^{\ell} \ye    &  \textrm{if $\ell' = 0$ and $\ell+\ell'$ is odd}  \\
       p_{+2}^{\ell} q_{+2}^{\ell'}\ve    &  \textrm{if $\ell' \ge 1$ and $\ell+\ell'$ is odd}
   \end{cases}
       \\[2mm]
   \sfb_{\ell,\ell'}
   & = &
   \begin{cases}
       p_{-1}^{\ell} \xe    &  \textrm{if $\ell' = 0$ and $\ell+\ell'$ is even}  \\
       p_{-1}^{\ell} q_{-1}^{\ell'}\ue    &  \textrm{if $\ell' \ge 1$ and $\ell+\ell'$ is even} \\
       p_{-2}^{\ell} \xo    &  \textrm{if $\ell' = 0$ and $\ell+\ell'$ is odd}  \\
       p_{-2}^{\ell} q_{-2}^{\ell'}\uo    &  \textrm{if $\ell' \ge 1$ and $\ell+\ell'$ is odd}
   \end{cases}
\end{subeqnarray}
Inserting these specializations into
Proposition~\ref{prop.cyclealt.Stype.final1.lambda=-1}, we obtain:

\begin{proposition}[$p,q$ S-fraction for cycle-alternating permutations, $\lambda=-1$]
   \label{prop.cyclealt.Stype.final1.lambda=-1.pq}
The ordinary generating function of the
polynomials~\reff{def.Qn.cycle-alternating.evenodd.pqgen} 
at $\lambda = -1$ has the S-type continued fraction
\begin{subeqnarray}
   \sum_{n=0}^\infty Q_n(\xe,\ye,\ue,\ve,\xo,\yo,\uo,\vo,p_{-1},p_{-2},p_{+1},p_{+2},q_{-1},q_{-2},q_{+1},q_{+2}, -1) \: t^n\nonumber\\
   \;=\;
	\cfrac{1}{1 + \cfrac{\xe \yo t}{1 +  \cfrac{(-p_{-2}\xo\!+\!q_{-2}\uo)(-p_{+2}\ye\!+\!q_{+2}\ve) t}{1 + \cfrac{(p_{-1}^2\xe\!+\!q_{-1}[2]_{-p_{-1},q_{-1}}\ue)(p_{+1}^2\yo\!+\!q_{+1}[2]_{-p_{+1},q_{+1}}\vo) t}{1 - \cdots}}}}
   \label{eq.thm.perm.Stype.cycle-alternating.evenodd.pqgen}
\end{subeqnarray}
with coefficients
\begin{subeqnarray}
   \alpha_{2k-1}  & = &  - (p_{-1}^{2k-2}\xe + q_{-1}[2k-2]_{-p_{-1},q_{-1}} \ue) \: (p_{+1}^{2k-2}\yo + q_{+1}[2k-2]_{-p_{+1},q_{+1}} \vo) \nonumber\\
	\mbox{}\\
   \alpha_{2k}    & = &  - (-p_{-2}^{2k-1}\xo + q_{-2}[2k-1]_{-p_{-2},q_{-2}} \uo) \: (-p_{+2}^{2k-1}\ye + q_{+2}[2k-1]_{-p_{+2},q_{+2}} \ve)\nonumber\\
	\mbox{}
 \label{def.weights.perm.Stype.cycle-alternating.evenodd.pqgen}
\end{subeqnarray}
\end{proposition}

Finally, denote by $Q_n(\xe,\ye,\ue,\ve,\xo,\yo,\uo,\vo,\lambda)$
the polynomial \reff{def.Qn.cycle-alternating.evenodd.pqgen} specialized to 
$p_{+1}=p_{+2}=p_{-1}=p_{-2}=q_{+1}=q_{+2}=q_{-1}=q_{-2} = 1$.
Setting $\lambda = -1$, we obtain:

\begin{proposition}[Simple S-fraction for cycle-alternating permutations, $\lambda=-1$]
   \label{prop.cyclealt.Stype.lambda=-1}
The ordinary generating function of the polynomials
$Q_n(\xe,\ye,\ue,\ve,\xo,\yo,\uo,\vo,-1)$ has the S-type continued fraction
\be
   \sum_{n=0}^\infty Q_n(\xe,\ye,\ue,\ve,\xo,\yo,\uo,\vo,-1) \: t^n
   \;=\;
   \cfrac{1}{1 + \cfrac{\xe \yo t}{1 +  \cfrac{(\xo\!-\!\uo)(\ye\!-\!\ve) t}{1 + \cfrac{\xe\yo t}{1 - \cdots}}}}
   \label{eq.thm.perm.Stype.cycle-alternating.evenodd}
\ee
with coefficients
\begin{subeqnarray}
   \alpha_{2k-1}  & = & - \xe  \yo \\[1mm]
   \alpha_{2k}    & = & - (\xo - \uo) \: (\ye - \ve)
 \label{def.weights.perm.Stype.cycle-alternating.evenodd}
\end{subeqnarray}
\end{proposition}

This proves the continued fraction
that was conjectured in \cite[eq.~(A.6)]{Deb-Sokal_cyclealt}.

\section{Results for D-permutations}   \label{sec.D-permutations}

We recall \cite{Lazar_20,Lazar_22,Lazar_23,Deb-Sokal_genocchi}
that a \textbfit{D-permutation} is a permutation of $[2n]$ satisfying\\
${2k-1 \le \sigma(2k-1)}$ and ${2k \ge \sigma(2k)}$ for all $k$;
D-permutations provide a combinatorial model
for the Genocchi and median Genocchi numbers.
We write $\dperm_{2n}$ for the set of D-permutations of $[2n]$.
Our running example 2,
\begin{eqnarray}
\sigma & = & 7\: 1\: 9\: 2\: 5\: 4\: 8\: 6\: 10\: 3\: 11\: 12\: 14\: 13\:
        \nonumber\\
       & = & (1,7,8,6,4,2)\,(3,9,10)\,(5)\,(11)\,(12)\,(13,14) \in \Sym_{14} \;,
\end{eqnarray}
is an example of a D-permutation.

We proceed in the same way as in the preceding section,
beginning with the ``master'' T-fraction
and then obtaining the others by specialization.

\subsection{Master T-fraction}

Following \cite[Section~3.4]{Deb-Sokal_genocchi},
we introduce six infinite families of indeterminates\\
$\bsfa = (\sfa_{\ell,\ell'})_{\ell,\ell' \ge 0}$,
$\bsfb = (\sfb_{\ell,\ell'})_{\ell,\ell' \ge 0}$,
$\bsfc = (\sfc_{\ell,\ell'})_{\ell,\ell' \ge 0}$,
$\bsfd = (\sfd_{\ell,\ell'})_{\ell,\ell' \ge 0}$,
$\bsfe = (\sfe_\ell)_{\ell \ge 0}$,
$\bsff = (\sff_\ell)_{\ell \ge 0}$
and then define the polynomials
\begin{eqnarray}
   & & \hspace*{-6mm}
   Q_n(\bsfa,\bsfb,\bsfc,\bsfd,\bsfe,\bsff,\lambda)
   \;=\;
       \nonumber \\[4mm]
   & &
   \sum_{\sigma \in \dperm_{2n}}
   \lambda^{\cyc(\sigma)}
   \prod\limits_{i \in \Cval(\sigma)}  \! \sfa_{\ucross(i,\sigma),\,\unest(i,\sigma)}
   \prod\limits_{i \in \Cpeak(\sigma)} \!\!  \sfb_{\lcross(i,\sigma),\,\lnest(i,\sigma)}
       \:\times
       \qquad\qquad
       \nonumber \\[1mm]
   & & \qquad\qquad\;\;\;\,
   \prod\limits_{i \in \Cdfall(\sigma)} \!\!  \sfc_{\lcross(i,\sigma),\,\lnest(i,\sigma)}
   \;
   \prod\limits_{i \in \Cdrise(\sigma)} \!\!  \sfd_{\ucross(i,\sigma),\,\unest(i,\sigma)}
       \:\times
       \qquad\qquad
       \nonumber \\[1mm]
   & & \qquad\qquad\;\;\,
   \prod\limits_{i \in \Evenfix(\sigma)} \!\!\! \sfe_{\psnest(i,\sigma)}
   \;
   \prod\limits_{i \in \Oddfix(\sigma)}  \!\!\! \sff_{\psnest(i,\sigma)}
   \;.
   \quad
 \label{def.Qn.firstmaster.Dperm}
\end{eqnarray}
(This is \cite[eq.~(3.30)]{Deb-Sokal_genocchi}
 with a factor $\lambda^{\cyc(\sigma)}$ included.)
Then the first master T-fraction for D-permutations
\cite[Theorem~3.11]{Deb-Sokal_genocchi}
handles the case $\lambda=1$:
it states that the ordinary generating function of the polynomials
$Q_n(\bsfa,\bsfb,\bsfc,\bsfd,\bsfe,1)$
has the T-type continued fraction
\be
   \sum_{n=0}^\infty Q_n(\bsfa,\bsfb,\bsfc,\bsfd,\bsfe,\bsff,1) \: t^n
   \;=\;
\Scale[0.95]{
   \cfrac{1}{1 - \sfe_0 \sff_0 t - \cfrac{\sfa_{00} \sfb_{00} t}{1 - \cfrac{(\sfc_{00} + \sfe_1)(\sfd_{00} + \sff_1) t}{1 - \cfrac{(\sfa_{01} + \sfa_{10})(\sfb_{01} + \sfb_{10}) t}{1 - \cfrac{(\sfc_{01} + \sfc_{10} + \sfe_2)(\sfd_{01} + \sfd_{10} + \sff_2)t}{1 - \cdots}}}}}
}
   \label{eq.thm.Tfrac.first.master}
\ee
with coefficients
\begin{subeqnarray}
   \alpha_{2k-1}
   & = &
   \biggl( \sum_{\ell=0}^{k-1} \sfa_{\ell,k-1-\ell} \biggr)
   \biggl( \sum_{\ell=0}^{k-1} \sfb_{\ell,k-1-\ell} \biggr)
        \\[2mm]
   \alpha_{2k}
   & = &
   \biggl( \sfe_k \,+\, \sum_{\ell=0}^{k-1} \sfc_{\ell,k-1-\ell} \biggr)
   \biggl( \sff_k \,+\, \sum_{\ell=0}^{k-1} \sfd_{\ell,k-1-\ell} \biggr)
        \\[2mm]
   \delta_1  & = &  \sfe_0 \sff_0 \\[2mm]
   \delta_n  & = &   0    \qquad\hbox{for $n \ge 2$}
 \label{def.weights.Tfrac.first.master}
\end{subeqnarray}
By Lemma~\ref{lemma1.1}, we obtain the case $\lambda=-1$
by inserting a factor $-1$ for each even or odd fixed point,
for each cycle peak (or alternatively, cycle valley),
and for each lower or upper crossing.
We therefore have:

\begin{proposition}[Master T-fraction for D-permutations, $\lambda = -1$]
   \label{thm.Tfrac.first.master}
The ordinary generating function of the polynomials
$Q_n(\bsfa,\bsfb,\bsfc,\bsfd,\bsfe,\bsff,-1)$
has the T-type continued fraction
\be
   \sum_{n=0}^\infty Q_n(\bsfa,\bsfb,\bsfc,\bsfd,\bsfe,\bsff,-1) \: t^n
   \;=\;
\Scale[0.95]{
   \cfrac{1}{1 - \sfe_0 \sff_0 t + \cfrac{\sfa_{00} \sfb_{00} t}{1 - \cfrac{(\sfc_{00} - \sfe_1)(\sfd_{00} - \sff_1) t}{1 + \cfrac{(\sfa_{01} - \sfa_{10})(\sfb_{01} - \sfb_{10}) t}{1 - \cfrac{(\sfc_{01} - \sfc_{10} - \sfe_2)(\sfd_{01} - \sfd_{10} - \sff_2)t}{1 - \cdots}}}}}
}
   \label{eq.thm.Tfrac.first.master.lambda=-1}
\ee
with coefficients
\begin{subeqnarray}
   \alpha_{2k-1}
   & = &
   -\, \biggl( \sum_{\ell=0}^{k-1} (-1)^\ell \, \sfa_{\ell,k-1-\ell} \biggr)
   \biggl( \sum_{\ell=0}^{k-1} (-1)^\ell \, \sfb_{\ell,k-1-\ell} \biggr)
        \\[2mm]
   \alpha_{2k}
   & = &
   \biggl( -\sfe_k \,+\, \sum_{\ell=0}^{k-1} (-1)^\ell \, \sfc_{\ell,k-1-\ell} \biggr)
   \biggl( -\sff_k \,+\, \sum_{\ell=0}^{k-1} (-1)^\ell \, \sfd_{\ell,k-1-\ell} \biggr)
        \\[2mm]
   \delta_1  & = &  \sfe_0 \sff_0 \\[2mm]
   \delta_n  & = &   0    \qquad\hbox{for $n \ge 2$}
 \label{def.weights.Tfrac.first.master.lambda=-1}
\end{subeqnarray}
\end{proposition}

\subsubsection{Running example 2}

We now write out the
monomial contributed by our running example~2
to the polynomial $Q_n(\bsfa,\bsfb,\bsfc,\bsfd,\bsfe,\bsff,\lambda)$
in \reff{def.Qn.firstmaster.Dperm}.
We have
$\sigma  =  7\: 1\: 9\: 2\: 5\: 4\: 8\: 6\: 10\: 3\: 11\: 12\: 14\: 13\\
        =  (1,7,8,6,4,2)\,(3,9,10)\,(5)\,(11)\,(12)\,(13,14) \in \dperm_{14}$,
which was depicted in Figure~\ref{fig.pictorial.2}.
Here $n=7$ and $\cyc(\sigma) = 6$.

The monomial contributed by $\sigma$ in~\reff{def.Qn.firstmaster.Dperm}
is almost the same as the monomial in~\reff{eq.monomial.example.2};
only the contribution of the fixed points is slightly different
because we treat even and odd fixed points separately.
Instead of
\be
\lambda^6\,
\sfa_{0,0}^2 \, \sfa_{1,0} \,
\sfb_{0,0}^2 \, \sfb_{0,1} \,
\sfc_{0,0} \, \sfc_{0,1} \, \sfc_{1,0} \,
\sfd_{0,0} \, \sfd_{0,1} \,
\sfe_{0}^2 \, \sfe_{2}
\ee
as in \reff{eq.monomial.example.2}, here the contribution is
\be
\lambda^6\,
\sfa_{0,0}^2 \, \sfa_{1,0} \,
\sfb_{0,0}^2 \, \sfb_{0,1} \,
\sfc_{0,0} \, \sfc_{0,1} \, \sfc_{1,0} \,
\sfd_{0,0} \, \sfd_{0,1} \,
\sfe_{0} \,
\sff_{0} \, \sff_{2}\;.
\ee


\subsection[\texorpdfstring{$p,q$}{p,q} T-fraction]{{\texorpdfstring{$\bm{p,q}$}{p,q}} T-fraction}

Consider now the polynomial
\begin{eqnarray}
   & &
   \hspace*{-14mm}
   P_n(x_1,x_2,y_1,y_2,u_1,u_2,v_1,v_2,\we,\wo,\ze,\zo,p_{-1},p_{-2},p_{+1},p_{+2},q_{-1},q_{-2},q_{+1},q_{+2},\se,\so,\lambda)
   \;=\;
   \hspace*{-1cm}
       \nonumber \\[4mm]
   & & \qquad\qquad
   \sum_{\sigma \in \dperm_{2n}}
   x_1^{\eareccpeak(\sigma)} x_2^{\eareccdfall(\sigma)}
   y_1^{\ereccval(\sigma)} y_2^{\ereccdrise(\sigma)}
   \:\times
       \qquad\qquad
       \nonumber \\[-1mm]
   & & \qquad\qquad\qquad\:
   u_1^{\nrcpeak(\sigma)} u_2^{\nrcdfall(\sigma)}
   v_1^{\nrcval(\sigma)} v_2^{\nrcdrise(\sigma)}
   \:\times
       \qquad\qquad
       \nonumber \\[3mm]
   & & \qquad\qquad\qquad\:
   \we^{\evennrfix(\sigma)} \wo^{\oddnrfix(\sigma)}
   \ze^{\evenrar(\sigma)} \zo^{\oddrar(\sigma)}
   \:\times
       \qquad\qquad
       \nonumber \\[3mm]
   & & \qquad\qquad\qquad\:
   p_{-1}^{\lcrosscpeak(\sigma)}
   p_{-2}^{\lcrosscdfall(\sigma)}
   p_{+1}^{\ucrosscval(\sigma)}
   p_{+2}^{\ucrosscdrise(\sigma)}
          \:\times
       \qquad\qquad
       \nonumber \\[3mm]
   & & \qquad\qquad\qquad\:
   q_{-1}^{\lnestcpeak(\sigma)}
   q_{-2}^{\lnestcdfall(\sigma)}
   q_{+1}^{\unestcval(\sigma)}
   q_{+2}^{\unestcdrise(\sigma)}
          \:\times
       \qquad\qquad
       \nonumber \\[3mm]
   & & \qquad\qquad\qquad\:
   \se^{\epsnest(\sigma)}
   \so^{\opsnest(\sigma)}
   \lambda^{\cyc(\sigma)}\;.
 \label{def.Pn.pq.Dperm}
\end{eqnarray}
(This is \cite[eq.~(3.22)]{Deb-Sokal_genocchi}
with a factor $\lambda^{\cyc(\sigma)}$ included.)
The various statistics have been defined in
\cite[Sections~2.7 and 2.8 and eq.~(3.22)]{Deb-Sokal_genocchi}.
This polynomial is obtained from \reff{def.Qn.firstmaster.Dperm}
by making the specializations \cite[eqs.~(6.40)-(6.45)]{Deb-Sokal_genocchi}
\begin{subeqnarray}
   \sfa_{\ell,\ell'}
   & = &
   p_{+1}^\ell q_{+1}^{\ell'}  \,\times\,
   \begin{cases}
      y_1  & \textrm{if $\ell' = 0$}   \\
      v_1  & \textrm{if $\ell' \ge 1$}
   \end{cases}
        \label{eq.proof.weights.sfa}  \\[2mm]
   \sfb_{\ell,\ell'}
   & = &
   p_{-1}^\ell  q_{-1}^{\ell'}  \,\times\,
   \begin{cases}
      x_1  & \textrm{if $\ell' = 0$}   \\
      u_1  & \textrm{if $\ell' \ge 1$}
   \end{cases}
        \\[2mm]
   \sfc_{\ell,\ell'}
   & = &
   p_{-2}^\ell q_{-2}^{\ell'}  \,\times\,
   \begin{cases}
      x_2  & \textrm{if $\ell' = 0$}   \\
      u_2  & \textrm{if $\ell' \ge 1$}
   \end{cases}
        \\[2mm]
   \sfd_{\ell,\ell'}
   & = &
   p_{+2}^\ell q_{+2}^{\ell'}  \,\times\,
   \begin{cases}
      y_2  & \textrm{if $\ell' = 0$}   \\
      v_2  & \textrm{if $\ell' \ge 1$}
   \end{cases}
        \\[2mm]
   \sfe_k
   & = &
   \begin{cases}
      \ze  & \textrm{if $k = 0$}   \\[1mm]
      \se^k \we  & \textrm{if $k \ge 1$}
   \end{cases}
        \\[2mm]
   \sff_k
   & = &
   \begin{cases}
      \zo  & \textrm{if $k = 0$}   \\[1mm]
      \so^k \wo  & \textrm{if $k \ge 1$}
   \end{cases}
        \label{eq.proof.weights.sff}
\end{subeqnarray}
Making these specializations in
Proposition~\ref{thm.Tfrac.first.master}
--- or equivalently, attaching a minus sign to the variables
$x_1, u_1, p_{+1}, p_{+2}, p_{-1}, p_{-2}, \we, \wo, \ze, \zo$
in \cite[Theorem~3.9]{Deb-Sokal_genocchi} ---
we obtain:


\begin{proposition}[$p,q$ T-fraction for D-permutations, $\lambda=-1$]
   \label{thm.Tfrac.first.pq.Dperm}
The ordinary generating function of the polynomials \reff{def.Pn.pq.Dperm}
at $\lambda = -1$ has the T-type continued fraction
\begin{eqnarray}
   & & \hspace*{-7mm}
\Scale[0.92]{
   \sum\limits_{n=0}^\infty
   P_n(x_1,x_2,y_1,y_2,u_1,u_2,v_1,v_2,\we,\wo,\ze,\zo,p_{-1},p_{-2},p_{+1},p_{+2},q_{-1},q_{-2},q_{+1},q_{+2},\se,\so,-1)  \: t^n
   \;=\;
}
       \nonumber \\[2mm]
   & & \hspace*{-3mm}
\Scale[0.87]{
   \cfrac{1}{1 - \ze \zo  \,t + \cfrac{x_1 y_1  \,t}{1 -  \cfrac{(x_2\!-\!\se\we)(y_2\!-\!\so\wo) \,t}{1 + \cfrac{(-p_{-1}x_1\!+\!q_{-1}u_1)(-p_{+1}y_1\!+\!q_{+1}v_1)  \,t}{1 - \cfrac{(-p_{-2}x_2\!+\!q_{-2}u_2\!-\!\se^2\we)(-p_{+2}y_2\!+\!q_{+2}v_2\!-\!\so^2\wo)  \,t}{1 + \cfrac{(p_{-1}^2 x_1\!+\! q_{-1} [2]_{-p_{-1},q_{-1}}u_1)(p_{+1}^2 y_1\!+\! q_{+1} [2]_{-p_{+1},q_{+1}}v_1)  \,t}{1 - \cfrac{(p_{-2}^2 x_2\!+\! q_{-2} [2]_{-p_{-2},q_{-2}} u_2\!-\!\se^3\we)(p_{+2}^2 y_2\!+\! q_{+2} [2]_{-p_{+2},q_{+2}}v_2\!-\!\so^3\wo)  \,t}{1 - \cdots}}}}}}}
}
       \nonumber \\[1mm]
   \label{eq.thm.Tfrac.first.pq.Dperm}
\end{eqnarray}
with coefficients
\begin{subeqnarray}
   \alpha_{2k-1} & = &
	\Scale[0.85]{-\bigl( (-p_{-1})^{k-1} x_1 + q_{-1} [k-1]_{-p_{-1},q_{-1}} u_1 \bigr) \:
	\bigl( (-p_{+1})^{k-1} y_1 + q_{+1} [k-1]_{-p_{+1},q_{+1}} v_1 \bigr)}
        \\[2mm]
   \alpha_{2k}   & = &
	\Scale[0.85]{\bigl( (-p_{-2})^{k-1} x_2 + q_{-2} [k-1]_{-p_{-2},q_{-2}} u_2 - \se^k \we \bigr) \:
	\bigl( (-p_{+2})^{k-1} y_2 + q_{+2} [k-1]_{-p_{+2},q_{+2}} v_2 - \so^k \wo \bigr)}
        \nonumber \\ \\
   \delta_1  & = &   \ze \zo   \\[1mm]
   \delta_n  & = &   0    \qquad\hbox{for $n \ge 2$}
   \label{eq.thm.Tfrac.first.weights.pq.Dperm}
\end{subeqnarray}
\end{proposition}

\subsubsection{Running example 2}

We now write out the monomial contributed by our running example~2
to the polynomial \reff{def.Pn.pq.Dperm} for $n=7$.
We have
$\sigma  =  7\: 1\: 9\: 2\: 5\: 4\: 8\: 6\: 10\: 3\: 11\: 12\: 14\: 13\\
        =  (1,7,8,6,4,2)\,(3,9,10)\,(5)\,(11)\,(12)\,(13,14) \in \dperm_{14}$,
which was depicted in Figure~\ref{fig.pictorial.2}.
Here $n=7$ and $\cyc(\sigma) = 6$.

The monomial contributed by $\sigma$ in~\reff{def.Pn.pq.Dperm}
is almost the same as the monomial in~\reff{eq.monomial.example.2.pq};
the contribution of the fixed points is slightly different
because we treat even and odd fixed points separately,
and because \reff{def.Qn.BIG.cyc.pq} distinguished fixed points by level
(subscripts on $w$), which we do not do here except to distinguish
level 0 (rar) from level $>0$ (nrfix).
Therefore, instead of
\be
\lambda^6 \:
x_1^2  \:
x_2^2 \:
y_1^3 \:
y_2 \:
u_1 \:
u_2 \:
v_2 \:
w_0^2 \:
w_2 \:
p_{+1} \:
p_{-2} \:
q_{+2} \:
q_{-1} \:
q_{-2} \:
s^2
\ee
as in \reff{eq.monomial.example.2.pq}, here the contribution is
\be
\lambda^6 \:
x_1^2 \:
x_2^2\:
y_1^3 \:
y_2 \:
u_1 \:
u_2\:
v_2 \:
\wo \:
\ze \: \zo\:
p_{+1} \:
p_{-2} \:
q_{+2} \:
q_{-1} \:
q_{-2} \:
\so^2\;.
\label{eq.monomial.example.2.pq.bis}
\ee

\subsection{Simple T-fraction}

Finally, denote by
$P_n(x_1,x_2,y_1,y_2,u_1,u_2,v_1,v_2,\we,\wo,\ze,\zo,\lambda)$
the polynomial \reff{def.Pn.pq.Dperm} specialized to 
$p_{+1}=p_{+2}=p_{-1}=p_{-2}=q_{+1}=q_{+2}=q_{-1}=q_{-2} = \se = \so = 1$.
This polynomial was introduced in \cite[eq.~(4.2)]{Deb-Sokal_genocchi}.
Making this same specialization in
Proposition~\ref{thm.Tfrac.first.pq.Dperm}
and using \reff{eq.n.minus.one.pq}, we obtain:


\begin{proposition}[Simple T-fraction for D-permutations, $\lambda=-1$]
   \label{thm.Tfrac.first.Dperm}
The ordinary generating function of the polynomials
$P_n(x_1,x_2,y_1,y_2,u_1,u_2,v_1,v_2,\we,\wo,\ze,\zo,-1)$
has the T-type continued fraction
\begin{eqnarray}
   & & \hspace*{-12mm}
   \sum_{n=0}^\infty
   P_n(x_1,x_2,y_1,y_2,u_1,u_2,v_1,v_2,\we,\wo,\ze,\zo,-1) \: t^n
   \;=\;
       \nonumber \\
   & &
   \cfrac{1}{1 - \ze \zo  \,t + \cfrac{x_1 y_1  \,t}{1 -  \cfrac{(x_2\!-\!\we)(y_2\!-\!\wo) \,t}{1 + \cfrac{(x_1\!-\!u_1)(y_1\!-\!v_1)  \,t}{1 - \cfrac{(x_2\!-\!u_2\!+\!\we)(y_2\!-\!v_2\!+\!\wo)  \,t}{1 + \cfrac{x_1 y_1  \,t}{1 - \cfrac{(x_2\!-\!\we)(y_2\!-\!\wo) \,t}{1 - \cdots}}}}}}}
       \nonumber \\[1mm]
   \label{eq.thm.Tfrac.first.Dperm}
\end{eqnarray}
with coefficients
\begin{subeqnarray}
   \alpha_{2k-1} & = &
      \begin{cases}
        - x_1 y_1                    & \textrm{if $k$ is odd}  \\
        - (x_1 - u_1) (y_1 - v_1)    & \textrm{if $k$ is even}
      \end{cases}
        \\[1mm]
   \alpha_{2k}   & = &
	\begin{cases}
	 (x_2-\we) (y_2-\wo)                    & \textrm{if $k$ is odd}  \\
         (x_2 - u_2+\we) (y_2 - v_2+\wo)    & \textrm{if $k$ is even}
      \end{cases}
        \\[1mm]
   \delta_1  & = &   \ze \zo   \\[1mm]
   \delta_n  & = &   0    \qquad\hbox{for $n \ge 2$}
   \label{eq.thm.Tfrac.first.weights.Dperm}
\end{subeqnarray}
\end{proposition}

\medskip

Finally, as a special case of Proposition~\ref{thm.Tfrac.first.Dperm},
we can obtain a J-fraction that was conjectured in
\cite[Appendix, case $\lambda=-1$]{Deb-Sokal_genocchi}.
It suffices to specialize the polynomials
$P_n(x_1,x_2,y_1,y_2,u_1,u_2,v_1,v_2,\we,\wo,\ze,\zo,\lambda)$
by setting $x_1 = x_2 = \ze = \zo = x$, ${y_1 = y_2 = y}$,
and $u_1 = u_2 = v_1 = v_2 = \we = \wo = 1$;
this yields the polynomials
\be
   P_n(x,y,\lambda)
   \;=\;
   \sum_{\sigma\in \dperm_{2n}}
       x^{\arec(\sigma)}y^{\erec(\sigma)} \lambda^{\cyc(\sigma)}
 \label{eq.Pn.xylam.Dperm}
\ee
that were introduced in \cite[eqs.~(4.1) and (A.1)]{Deb-Sokal_genocchi}.
Inserting this specialization in Proposition~\ref{thm.Tfrac.first.Dperm}
gives, for $\lambda = -1$, a T-fraction with coefficients
\begin{subeqnarray}
   \alpha_{2k-1} & = &
      \begin{cases}
        - x y                    & \textrm{if $k$ is odd}  \\
        - (x -1) (y - 1)    & \textrm{if $k$ is even}
      \end{cases}
        \\[1mm]
   \alpha_{2k}   & = &
        \begin{cases}
         (x-1) (y-1)                    & \textrm{if $k$ is odd}  \\
         x  y    & \textrm{if $k$ is even}
      \end{cases}
        \\[1mm]
   \delta_1  & = &   x^2   \\[1mm]
   \delta_n  & = &   0    \qquad\hbox{for $n \ge 2$}
\end{subeqnarray}
Using the even contraction for T-fractions with
$\delta_2 = \delta_4 = \delta_6 = \ldots = 0$
\cite[Proposition~2.1]{Deb-Sokal_genocchi},
we can rewrite this as a J-fraction:

\begin{corollary}
        \label{cor.Dperm}
The ordinary generating function of the polynomials \reff{eq.Pn.xylam.Dperm}
has the J-type continued fraction
\be
   \sum_{n=0}^\infty
   P_n(x,y,-1)
   \;=\;
	\cfrac{1}{1 - x(x\!-\!y)  \,t + \cfrac{ xy\:(x-1)(y-1)  \,t^2}{1 +  \cfrac{xy\:(x-1)(y-1)  \,t^2}{1 + \cfrac{xy\:(x-1)(y-1)  \,t^2}{1 + \cdots}}}}
\ee
with coefficients
\begin{subeqnarray}
   \gamma_{0} & = & x(x-y)  \\[1mm]
   \gamma_{n} & = & 0 \quad \hbox{for $n\geq 1$}\\[1mm]
   \beta_{n}   & = & -xy(x-1)(y-1)
\end{subeqnarray}
\end{corollary}

This J-fraction was conjectured in
\cite[Appendix, case $\lambda=-1$]{Deb-Sokal_genocchi}.

\section*{Acknowledgments}

One of us (B.D.)\ is partially supported by the DIMERS project ANR-18-CE40-0033.
He also wishes to thank Jakob Stein for helpful discussions
concerning the topology of the plane.

\end{document}